\documentclass[11pt]{llncs}
\usepackage[latin1]{inputenc}
\usepackage[T1]{fontenc}

\usepackage{threeparttable}
\usepackage{graphicx}

\usepackage{pslatex,graphicx,xspace}
\usepackage{amsmath,amssymb,amsfonts}
\usepackage{url}
\date{}
\let\oldtitle\title\let\oldauthor\author
\renewcommand{\title}[1]{\oldtitle{\bfseries #1}}
\renewcommand{\author}[1]{\oldauthor{\bfseries #1}}
\makeatletter
\renewenvironment{abstract}%
  {\vspace{\baselineskip}\hrule\par\vspace{.3\baselineskip}\noindent\textbf{\abstractname} --- }
  {\par\vspace{.3\baselineskip}\hrule}
\makeatother

\newcommand{\R}{\mathbb{R}}
\newcommand{\frosen}{\ensuremath{f_{\mathrm{Rosen}}}}

\newcommand{\felli}{\ensuremath{f_{\mathrm{elli}}}}
\newcommand{\sqrtsqrtfelli}{\ensuremath{f_{\mathrm{elli}}^{\mathrm{1/4}}}}

\newcommand{\ftarget}{\ensuremath{f_\mathrm{target}}}
\newcommand{\basis}{\ensuremath{\vec{\vec{B}}}}

\newcommand{\comm}[1]{}
\renewcommand{\vec}[1]{\mathchoice{\mbox{\boldmath$\displaystyle#1$}}
  {\mbox{\boldmath$\textstyle#1$}} {\mbox{\boldmath$\scriptstyle#1$}}
  {\mbox{\boldmath$\scriptscriptstyle#1$}}}

\newcommand{\x}{\ensuremath{\vec{x}}}
\newcommand{\y}{\ensuremath{\vec{y}}}

\newcommand{\ie}{{\textit{i.e.}}}

\newcommand{\nref}[1]{\nolinebreak[4]\hspace{0.2845em plus0.05em minus0.17em}\nolinebreak[4]\ref{#1}}

\graphicspath{{figures/}}

\title{Experimental Comparisons of Derivative Free Optimization Algorithms$^1$}

\author{A.~Auger$^{\dag,\ddag}$, N.~Hansen$^{\dag,\ddag}$, J. M. Perez Zerpa$^{\dag}$, 
 R.~Ros$^{\dag}$, M.~Schoenauer$^{\dag,\ddag}$}
\institute{$^{\dag}$ TAO Projetct-Team, INRIA Saclay -- Île-de-France\\
LRI, Bat 490 Univ. Paris-Sud 91405 Orsay Cedex France\\
$^\ddag$ Microsoft Research-INRIA Joint Centre\\
28 rue Jean Rostand, 91893 Orsay Cedex, France
}

\begin{document}
\maketitle\thispagestyle{empty}
%

\begin{abstract}
In this paper, the performances of the quasi-Newton BFGS algorithm, the NEWUOA derivative free optimizer, the Covariance Matrix Adaptation Evolution Strategy (CMA-ES), the Differential Evolution (DE) algorithm and Particle Swarm Optimizers (PSO) are compared experimentally on benchmark functions reflecting important challenges encountered in real-world optimization problems. Dependence of the performances in the conditioning of the problem and rotational invariance of the algorithms are in particular investigated.

\end{abstract}


\section{Introduction}

Continuous Optimization Problems (COPs) aim at finding the global optimum (or optima) of a real-valued function (aka {\em objective} function) defined over a (subset of) a real vector space. COPs commonly appear in everyday's life of many scientists, engineers and researchers from various disciplines, from physics to mechanical, electrical and chemical engineering to biology. Problems such as model calibration, process control, design of parameterized parts are routinely modeled as COPs.
Furthermore, in many cases, very little is known about the objective function. In the worst case, it is only possible to retrieve objective function values for given inputs, and in particular the user has no information about derivatives, or even about some weaker characteristics of the objective function (e.g. monotonicity, roughness, \ldots). This is the case, for instance, when the objective function is the output of huge computer programs ensuing from several years of development, or when experimental processes need to be run in order to compute objective function values. Such problems amount to what is called {\em Black-Box Optimization} (BBO).\footnotetext[1]{Invited Paper at the 8$^{th}$ International Symposium on Experimental Algorithms, June 3-6, 2009, Dortmund, Germany}

Because BBO is a frequent situation, many optimization methods (aka {\em search algorithms}) have been proposed to tackle BBO problems, that can be grossly classified in two classes: (i) deterministic methods include classical derivative-based algorithms, in which the derivative is numerically computed by finite differences, and more recent Derivative Free Optimization (DFO) algorithms \cite{Conn:1997}, like pattern search \cite{Torczon:97} and trust region methods \cite{Powell:2006}; (ii) stochastic methods rely on random variables sampling to better explore the search space, and include recently introduced bio-inspired algorithms (see Section \ref{algorithms}).

However, the practitioner facing a BBO problem has to choose among those methods, and there exists no theoretical solid ground where he can stand to perform this choice, first because he does not know much about his objective function, but also because all theoretical results either make simplifying hypotheses that are not valid for real-world problems, or give results that do not yield any practical outcome. 
Moreover, most of BBO methods require some parameter tuning, and here again very little help is available for the practitioner, who is left with a blind and time-consuming test-and-trial approach.

In such context, this paper proposes an experimental perspective on BBO algorithms comparisons. Rigorous procedures to compare the results of different BBO algorithms have been proposed \cite{CEC05}, taking into account the stochastic nature of many of them, and giving fair chances to each one of them. However, a critical issue in such experiments is that of the benchmark suite. And because no set of real-world problems can be guaranteed to cover all possible cases of difficult COPs, the approach that has been chosen here is to build artificial test functions with some precise characteristics that are known to be possible sources of difficulty for optimization (e.g. ill-conditioning, non-separability, non-convexity, ruggedness, \ldots). Such experimental results could then be cautiously generalized, leaving only a few good-performing algorithms in each specific context.

Of course, in real-life BBO situations, it is assumed that nothing is known about the objective function. However, the user sometimes has some partial information (e.g. because his problem is known to be similar to other better-known problems) that might lead him to decide for a BBO method that is (experimentally) known to perform well, 'in vitro', in his precise situation. But on the other hand, assuming absolutely nothing is known in advance about the objective function, running the champion algorithms as identified in perfectly controlled environment might give him some information about his function (e.g. if numerical gradient-based algorithms perform 100 times better than all other methods, his problem is probably very similar to a quadratic problem). This paper is a first step in aiming such 'in vitro' results.

Next, in Section \ref{difficulty}, some characteristics of the objective function are surveyed that are known to make the corresponding BBO problem hard. Section \ref{algorithms} introduces the algorithms that will be compared here. Section \ref{functions} then introduces the test bench that illustrates the different difficulties highlighted in Section \ref{difficulty}, as well as the experimental conditions of the comparisons. The results are presented and discussed in Section \ref{results}, and the paper ends with some conclusions in Section \ref{conclusions}.

\section{What makes a search problem difficult?}
\label{difficulty}
In this section, we discuss problem characteristics that are especially challenging for search algorithms.

\subsection{Ill-conditioning}


The conditioning of a problem can be defined as the range (over a level set) of the maximum improvement of objective function value in a ball of small radius centered on a given level set. In the 
case of convex quadratic functions ($f(x)=\frac12 x^{T} H x$ where $H$ is a symmetric definite 
matrix), the conditioning can be exactly defined as the condition number of the Hessian matrix $H$, i.e., the ratio between the largest and smallest eigenvalue. Since level sets associated to a 
convex quadratic function are ellipsoids, the condition number corresponds to the squared ratio between the largest and shortest axis lengths of the ellipsoid. 

Problems are typically considered as ill-conditioned if the conditioning is larger than $10^{5}$. In practice we have seen problems with conditioning as large as $10^{10}$. In this paper we will quantitatively assess the performance dependency on the conditioning of the objective function.

\subsection{Non-separability}

An objective function $f(x_{1},\ldots, x_{n})$ is separable if the optimal value for any variable $x_{i}$ can be obtained by optimizing $f(\widetilde{x}_{1},\ldots, \widetilde{x}_{i-1} ,x_{i}, \widetilde{x}_{i+1},\ldots, \widetilde{x}_{n})$ for any fixed choice of the variables $\widetilde{x}_{1},\ldots,\widetilde{x}_{i-1},\widetilde{x}_{i+1},\ldots,\widetilde{x}_{n}$. Consequently optimizing an $n$-dimensional separable objective function reduces to optimizing $n$ one-dimensional functions.

Functions that are additively decomposable, i.e., that can be written as $f(\x) = \sum_{i=1}^{n} f_{i}(x_{i})$ are separable. One way to render a separable test function non-sepa\-rable is to rotate first the vector $\x$, which can be achieved by multiplying $\x$ by  an orthogonal matrix $\basis$: if $\x \mapsto f(\x)$ is separable, the function $\x \mapsto f(\basis \x)$ might be non-separable for all non-identity orthogonal matrices $\basis$. In this paper we will investigate separable and non-separable problems. 

\subsection{Non-convexity}
\label{convexity}
Some BBO methods implicitly assume or exploit convexity of the objective function. Composing a convex function $f\in \R$  to the left with a monotonous transformation $g: \R \rightarrow \R$ can result in a non-convex function, for instance the one-dimensional convex function $f(x)=x^{2}$ composed with the monotonous function $g(.) = |.|^{1/4}$ becomes the non-convex function $\sqrt{|.|}$. In this paper we will assess performance dependency on convexity.

\section{Algorithms tested}
\label{algorithms}
This section introduces the different algorithms that will be compared in this paper. They have been chosen because they are considered to be the champions in their category, both in the deterministic optimization world (BFGS and NEWUOA) and in the stochastic bio-inspired world (CMA-ES, DE and PSO). They will also be a priori discussed here with respect to the difficulties of continuous optimization problems highlighted in previous Section \ref{difficulty}.

\subsection{The algorithms}
\subsubsection{BFGS}
is a well-known quasi-Newton (i.e. gradient-based) method: from the current point, it computes a 'descent direction' using an approximation of the inverse of the Hessian matrix of the objective function applied to its gradient, and performs a line-search (1D optimization) along this direction. It then updates the approximate inverse Hessian. BFGS method is a local method: it has a proven convergence to a stationary point\ldots provided the starting point is close enough from the solution, and the objective function is regular. The Matlab$^{\mbox{\textregistered}}$ version of BFGS (Matlab function {\tt fminunc}) will be used here, because it is blindly used by many scientists facing optimization problems. Default parameters were used except for stopping criteria: the algorithms stops if the function value improvement in one iteration is less than $10^{-25}$.

In BBO context, the gradients have to be computed numerically (an option in Matlab BFGS), which might be a source of possible numerical problems. 

\subsubsection{NEWUOA}
(NEW Unconstrained Optimization Algorithm) has been proposed by Powell \cite{Powell:2006}: it is a DFO algorithm using the trust region paradigm. The trust region is a ball, centered on the current best point. NEWUOA computes a quadratic interpolation of the objective function within the current trust region, based on known values of the objective, and then performs a truncated conjugate gradient minimization of the interpolated model in the trust region. It then updates either the current best point or the radius of the trust region, based on the a posteriori interpolation error, and some thresholds on the trust region size. Here, the implementation by Matthieu Guibert posted at \url{http://www.inrialpes.fr/bipop/people/guilbert/newuoa/newuoa.html} has been used.

An important parameter of NEWUOA is the quadratic model to use for the interpolation, or, equivalently, the number of points that are necessary to compute the interpolation. As recommended by Powell \cite{Powell:2006}, $2n+1$ points have been used here ($n$ is the dimension of the search space). Other critical parameters are the initial and final radii of the trust region: the initial radius governs the granularity of the objective function that the algorithm will 'see' and the final radius tunes the amount of local search that will performed. Here the initial and final values 100 and $10^{-15}$ were used, after some preliminary experiments.

\subsubsection{CMA-ES}
is an Evolution Strategy (ES) \cite{rechenberg:73,Schwefel:1995} algorithm: from a set of 'parents' (potential solutions), 'offspring' are created by sampling Gaussian distributions, and the best of the offspring (according to the objective function values) become the next parents. The art of Evolution Strategies lies in the way the parameters of the Gaussian distributions are updated: the Covariance Matrix Adaptation \cite{Hansen:2001} uses the path that has been followed by evolution so far to (i) adapt the step-size, a scaling parameter that tunes the granularity of the search, comparing the actual path length to that of a random walk; and (ii) modify the covariance matrix of the multivariate Gaussian distribution by modifying its eigenvectors in order to increase the likelihood of recent beneficial moves. A single Gaussian distribution is maintained, its mean being a linear combination of the parents. Besides the population size, CMA-ES is parameter-free. The population size has been set to its default value $4 + \lfloor3 \log(n)\rfloor$, but it needs to be increased in order to tackle highly rugged search landscapes. The initial step-size has been set to a third of the parameters' range.
The version used in this paper (Scilab 0.92) implements weighted recombination and rank-$\mu$ update \cite{Hansen:2006} (version 0.99 is available at \url{http://www.lri.fr/~hansen/cmaes_inmatlab.html}).

\subsubsection{PSO}
(Particle Swarm Optimization) \cite{kennedy1995pso} is a bio-inspired algorithm that recently raised a lot of interest, thanks to several published good results, and the simplicity of its implementation. The biological paradigm is that of a swarm of particles that 'fly' over the objective landscape, exchanging information about the best locations (i.e. potential solutions) they have seen. More precisely, each particle updates its velocity, stochastically twisting it toward the direction of the best positions so far visited by (i) itself and (ii) the whole swarm; it then updates its position according to its velocity and computes the new value of the objective function.

A Scilab transcription of the Standard PSO 2006, that is still available on the main page of {\em PSO Central} \url{http://www.particleswarm.info/}, was used here, with default settings. 

\subsubsection{Differential Evolution}
(DE \cite{Storn:1997}) borrows from Evolutionary Algorithms (EAs) a population of potential solutions that evolves subject to objective-function based selection. However, the main operator used to generate new solutions, that somehow replaces mutation, is specific to DE (and the source for its name): the difference between two points in the population is added to a third one. Uniform crossover is used with some probability. The implementation posted by the original authors at \url{http://www.icsi.berkeley.edu/~storn/code.html} was used here. However, the authors themselves confess, in their guidance to DE parameter tuning, that the results might be very dependent on the parameters. They propose in the code 6 possible settings, and extensive experiments ($3\times288$ trials) on a moderately ill-conditioned problem lead us to consider the ``{\em DE/local-to-best/1/bin}'' strategy, where a single difference vector, computed between a random point and the best point in the population, is used to generate the new points. In those preliminary experiments, the use of crossover seemed to have little beneficial impact on the results, so no crossover was used, thus making DE rotationally invariant (see below).
Those preliminary experiments also resulted in values of the other parameters of DE: the population size was set to the recommended value of $10n$, a weighting factor to $F=0.8$.

\subsection{Invariances}
\label{invariances}
Some a priori comparisons can be made about those algorithms, related to the notion of {\em invariance}. Indeed, invariances add to the robustness of an algorithm: functions belonging to the same equivalence class with respect to some invariance property will look exactly the same for an algorithm that is invariant under the transformation defining this equivalence class.

Two sets of invariance properties are distinguished, whether they regard transformations of the objective function value or transformations of the search space.
First, all comparison-based algorithms are invariant under monotonous transformations of the objective function, as comparisons are unaltered if the objective function $f$ is replaced with some $g \circ f$ for some monotonous function $g$. All bio-inspired algorithms used in this paper are comparison-based, while the BFGS and NEWUAO are not (see Section \ref{convexity}).

Regarding transformations of the search space, all algorithms are trivially invariant under translation of the coordinate system. But let us consider some orthogonal rotations: BFGS is coordinate-dependent due to the computation of numerical gradients. NEWUOA is invariant under rotation when considering the complete quadratic model, i.e.\ built with $\frac{1}{2}(n+1)(n+2)$ points. This variant is however often more costly compared to the $2n+1$ one -- but the latter is not invariant under rotation. The rotational invariance of CMA-ES is built-in, while that of DE depends whether or not crossover is used -- as crossover relies on the coordinate system. This was one reason for omitting crossover here. Finally, PSO is (usually) not invariant under rotations, as all computations are done coordinate by coordinate \cite{hansen2008pso,wilke2007b:cla}.

\section{Test functions and experimental setup}
\label{functions}
\subsection{Test functions}
The benchmark functions tested are given in
Table\nref{testfuncs}.
\begin{table*}[tbh] 
  \caption[testfuncs]{\label{testfuncs} Test functions with coordinate-wise
  initialization intervals and target function value, where
  $\y:=\basis\x$ implements
  an angle-preserving, linear transformation, \ie\ \basis\ is
  orthogonal. }
\begin{center}
\begin{tabular}{lccc} 
Function                  & $\alpha$     & Initialization     & \ftarget
   \\ 
\hline 
$\felli(\x)   =
   \sum_{i=1}^{n}\alpha^{\frac{i-1}{n-1}}y_i^2 \comm{= \y'\left(
   \begin{array}{cccccc}
        \alpha^\frac{0}{n-1} & \\
        & \alpha^\frac{1}{n-1} &  \\
         &                      &\ddots & \\
          &                       &       & \alpha^\frac{n-1}{n-1} \\
   \end{array}
   \right)\y}$ & 
   $[1,10^{10}]$  &
   $[-20 , 80]^n$  &
   $10^{-9}$
   \\ 
$\frosen(\x)   = \sum_{i=1}^{n-1}\big( \alpha \, 
                   (y_i^2-y_{i+1})^2+ (y_i-1)^2\big)$ & 
                     $[1,10^8]$ &
                     $[-20 , 80]^n$      &
   $10^{-9}$
\\
$\sqrtsqrtfelli(\x)   =
   \left( \sum_{i=1}^{n}\alpha^{\frac{i-1}{n-1}}y_i^2 \right)^{1/4} \comm{= \y'\left(
   \begin{array}{cccccc}
        \alpha^\frac{0}{n-1} & \\
        & \alpha^\frac{1}{n-1} &  \\
         &                      &\ddots & \\
          &                       &       & \alpha^\frac{n-1}{n-1} \\
   \end{array}
   \right)\y}$ & 
   $[1,10^{10}]$  &
   $[-20 , 80]^n$  &
   $10^{-9}$
\end{tabular} 
\end{center}
\end{table*}
The functions are tested in their original axis-parallel version
(i.e. $\basis$ is the identity and $\y=\x$), and in rotated versions,
where $\y=\basis\x$. The orthogonal
matrix \basis\ is chosen such that each column is uniformly
distributed on the unit hypersphere surface \cite{Hansen:2001}, fixed
for each run.

 The ellipsoid function $\felli$ is a convex-quadratic function where the parameter $\alpha$
 is the condition number of the Hessian matrix that is varied
 between 1 and $10^{10}$ in our experiments. If $\alpha=1$
 the ellipsoid is the isotropic separable sphere function. The function $\sqrtsqrtfelli$ has the same contour lines (level sets) as $\felli$, however it is neither quadratic nor convex. For $\alpha\neq1$, the functions $\felli$ and $\sqrtsqrtfelli$ are separable if and only if
 $\basis=\vec{\vec{I}}$. 

The Rosenbrock function $\frosen$ is non-separable, has
 its global minimum at $\x=[1,1,\dots,1]$ and, for large enough
 $\alpha$ and $n$, has one local minimum close to $\x=[-1,1,\dots,1]$, see
 also \cite{shang2006ane}. The contour lines of the Rosenbrock function show a bent ridge that guides to the global optimum (the Rosenbrock is sometimes called banana function) and the parameter $\alpha$ controls the width of the ridge. In the classical
 Rosenbrock function $\alpha$ equals $100$. For smaller $\alpha$ the
 ridge becomes wider and the function becomes less difficult to
 solve. We vary $\alpha$ between one and $10^8$.

\subsection{Experimental Setup}

For each algorithm tested we conduct $21$ independent trials of up to $10^7$ function evaluations. If, for BFGS, no success was encountered, the number of trials was extended to 1001. 

We quantify the performance of the algorithms using the success performance $SP1$ used in \cite{Hansen:2004}, analyzed in \cite{Auger:2005b},
and also denoted as Q-measure in \cite{feoktistov:2006des}. The $SP1$  equals the average number of function evaluations for successful runs divided by the ratio of successful experiments, where an experiment is successful if the $\ftarget$ is reached before $10^{7}$ function evaluations are exceeded. The $SP1$ is an estimator of the expected number of function evaluations to reach $\ftarget$ if the algorithm is restarted until a success (supposing infinite time horizon) and assuming that the expected number of function evaluations for unsuccessful runs equals the expected number of evaluations for successful runs.

\section{Results}
\label{results}
\newcommand{\NEW}{NEWUOA}
\newcommand{\SElli}{Ellipsoid$^{1/4}$}
\newcommand{\includegraphtrim}[1]{%
    \includegraphics[width=1.15\textwidth,trim=0 0 0 1.5cm,clip]{figures/#1}}
Results are shown for dimension 20. Results for 10 and 40D reveal similar tendencies and are displayed in Appendix~\ref{app:results}. 

\paragraph{Ellipsoid functions: dependencies} 
%
\begin{figure}[tb]
\begin{minipage}{0.49\textwidth}\centering
Separable Ellipsoid Function \\ \includegraphtrim{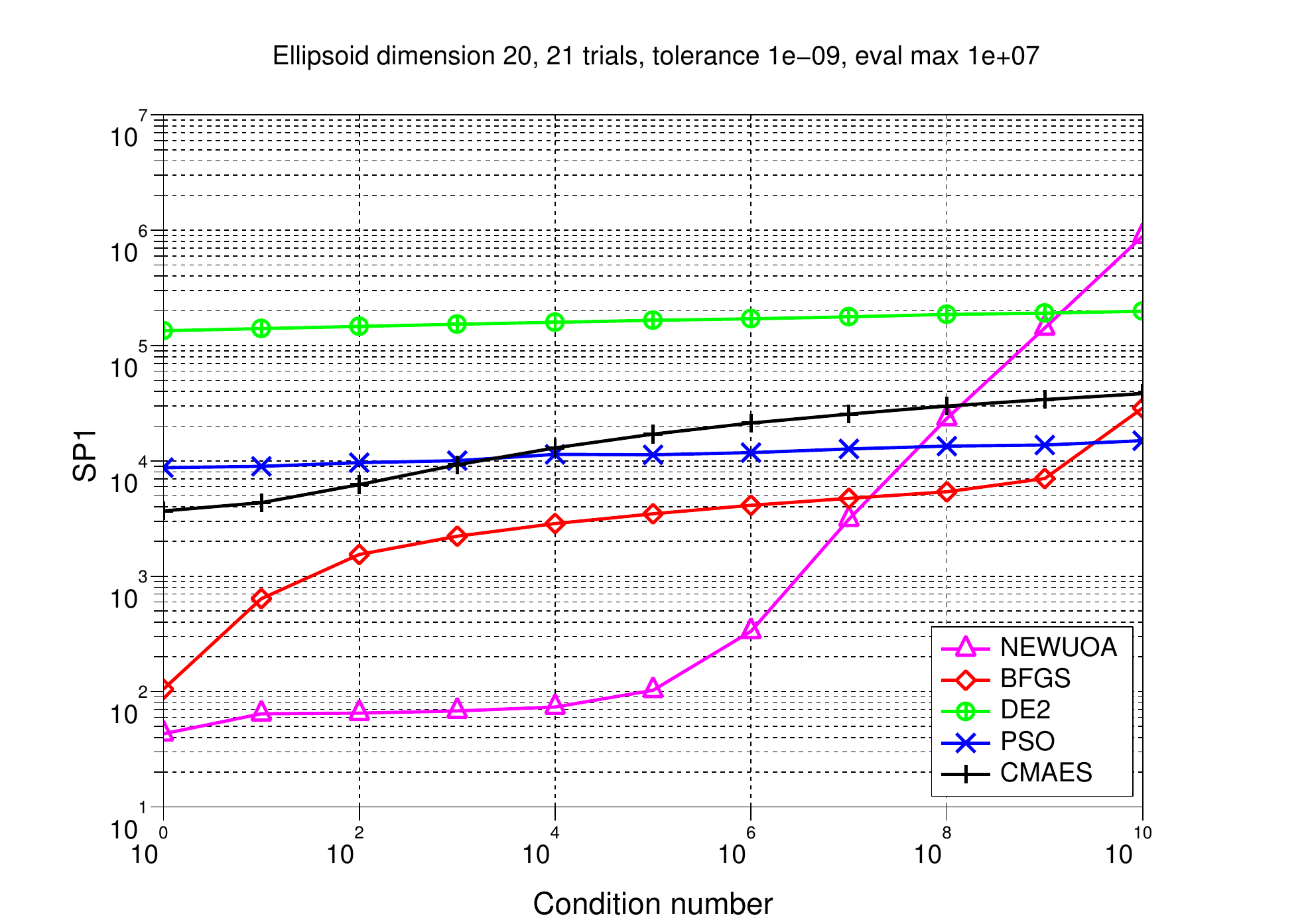}\\[1ex]
Separable \SElli\ Function \\ \includegraphtrim{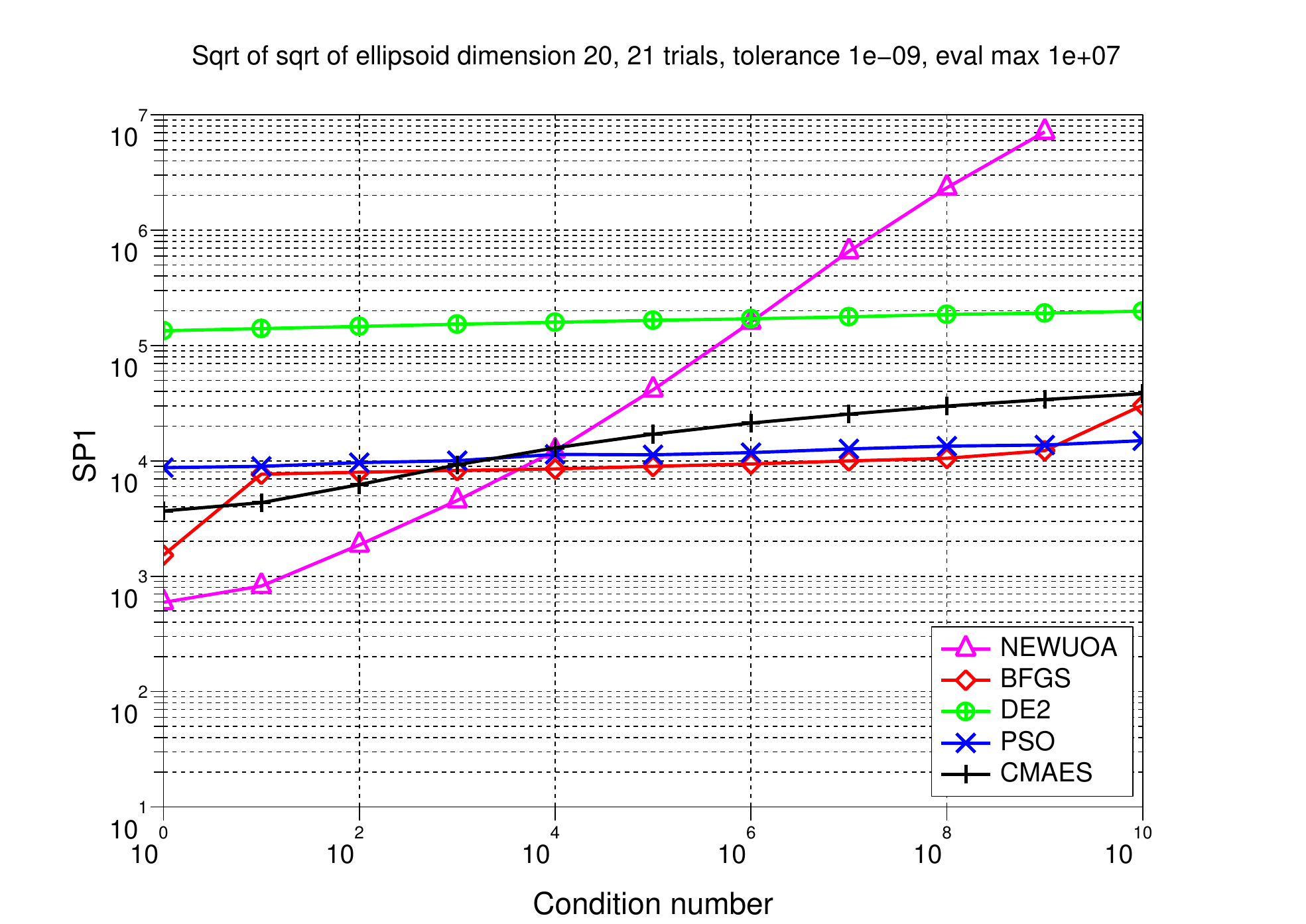}
\end{minipage}
\begin{minipage}{0.49\textwidth}\centering
Rotated Ellipsoid Function \\ \includegraphtrim{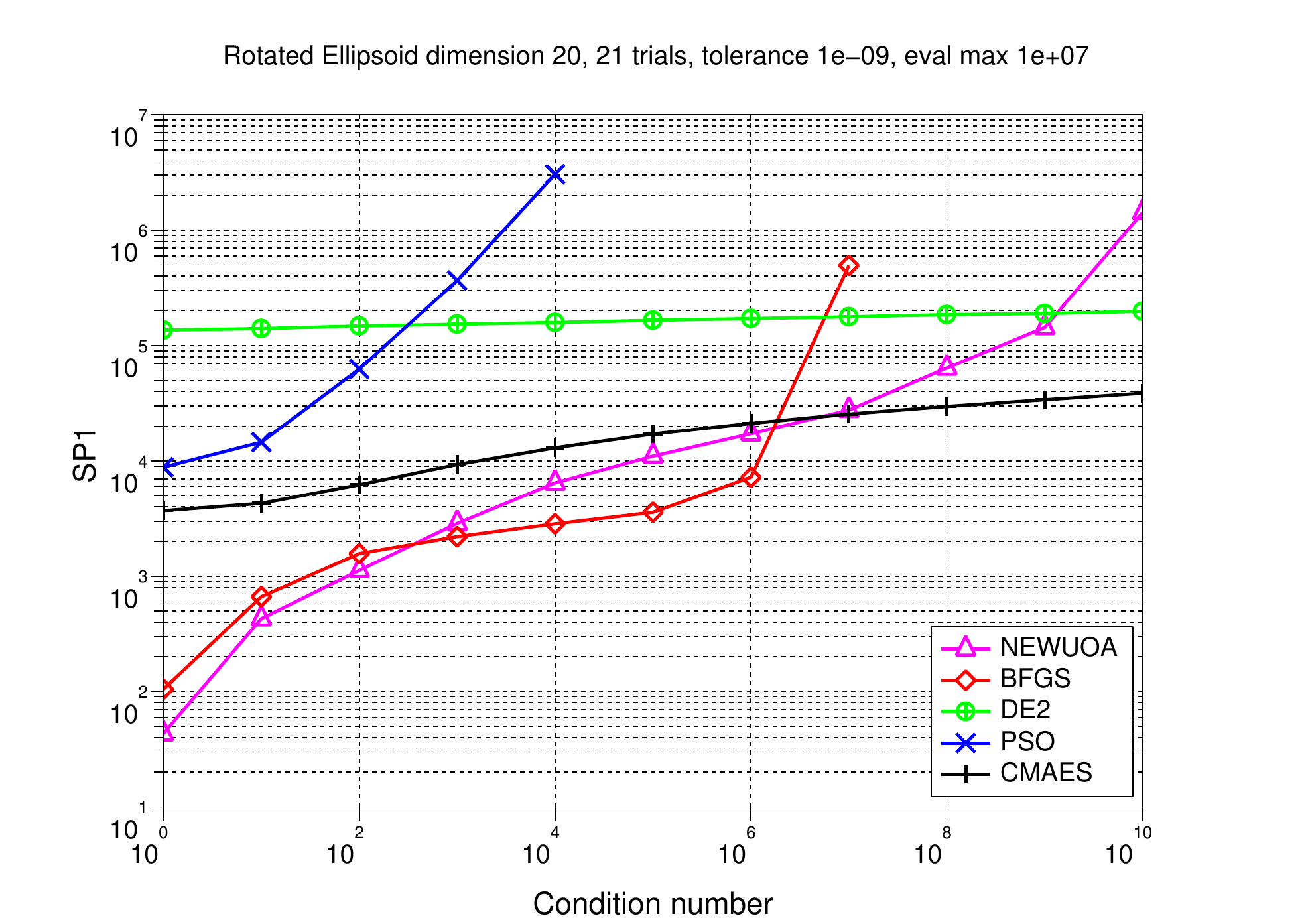}\\[1ex]
Rotated \SElli\ Function \\ \includegraphtrim{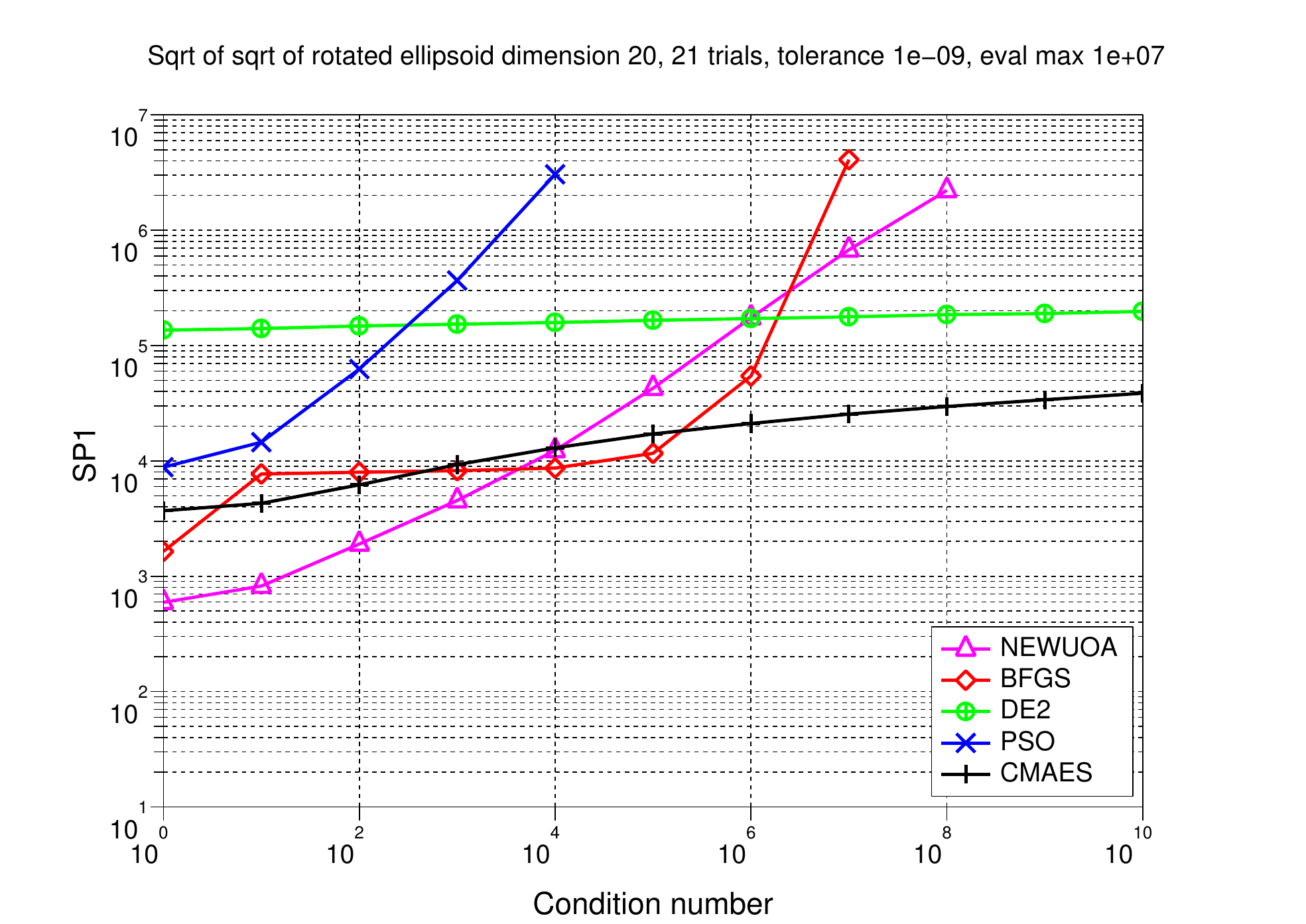}
\end{minipage}
\caption[elli]{\label{fig:ellitwenty}All ellipsoidal functions in 20D. Shown is $SP1$ (the expected running time  or number of function evaluations to reach the target function value) versus condition number. }
\end{figure}
%
Figure~\ref{fig:ellitwenty} shows $SP1$ (search costs, expected running time in number of function evaluations) versus condition number on all ellipsoidal functions. A remarkable dependency on the condition number can be observed in most cases. The two exceptions are PSO on the separable functions and DE. In the other cases the performance declines by at least a factor of ten for very ill-conditioned problems as for CMA-ES. The overall strongest performance decline is shown by PSO on the rotated functions. \NEW\ shows in general a comparatively strong decline, while BFGS is only infeasible for high condition numbers in the rotated case, reporting some numerical problems. The decline of CMA-ES is moderate. 

For CMA-ES and DE the results are (virtually) independent of the given ellipsoidal functions, where CMA-ES is consistently between five and forty times faster than DE. For PSO the results are identical on Ellipsoid and \SElli, while the performance decline under rotation (left versus right figures) is very pronounced. PSO performs well only on separable or very well-conditioned functions. A similar strong decline under rotation can be observed for \NEW\ on the Ellipsoid function for moderate condition numbers. BFGS, on the other hand, shows a strong rotational dependency on both functions only for large condition numbers $\ge10^6$. 

Switching from Ellipsoid (above) to \SElli\ (below) only effects BFGS and \NEW. BFGS becomes roughly five to ten times slower. A similar effect can be seen for \NEW\ on the rotated function. On the separable Ellipsoid function the effect is more pronounced, because \NEW\ performs exceptionally well on the separable Ellipsoid function.  

\paragraph{Ellipsoid functions: comparison} On the separable Ellipsoid function up to a condition number of $10^6$ \NEW\ clearly outperforms all other algorithms. Also BFGS performs still better than PSO and CMA-ES while DE performs worst. On the separable \SElli\ function BFGS, CMA-ES and PSO perform similar. \NEW\ is faster for low condition numbers and slower for large ones. For condition number larger than $10^6$, \NEW\ becomes even worse than DE.  

On the rotated functions, the performance of PSO declines fast with increasing condition number. For numbers larger than $10^3$, PSO is remarkably outperformed by all other algorithms.  On the rotated Ellipsoid function for moderate condition numbers BFGS and \NEW\ perform best and outperform CMA-ES by a factor of five, somewhat more for low condition numbers, and less for larger condition numbers, while PSO and DE are much worse. For large condition numbers CMA-ES becomes superior and DE is within a factor of ten of the best performance. 

On the rotated \SElli\ BFGS and CMA-ES perform similar up to condition of $10^6$. \NEW\ performs somewhat better for lower condition numbers up to $10^4$. For larger condition numbers BFGS and \NEW\ decline and CMA-ES performs best. 

\paragraph{Rosenbrock function}
\newcommand{\alp}{\ensuremath{\alpha}}
%
\begin{figure}[t]
\begin{minipage}{0.49\textwidth}\centering
Rosenbrock Function\\
\includegraphtrim{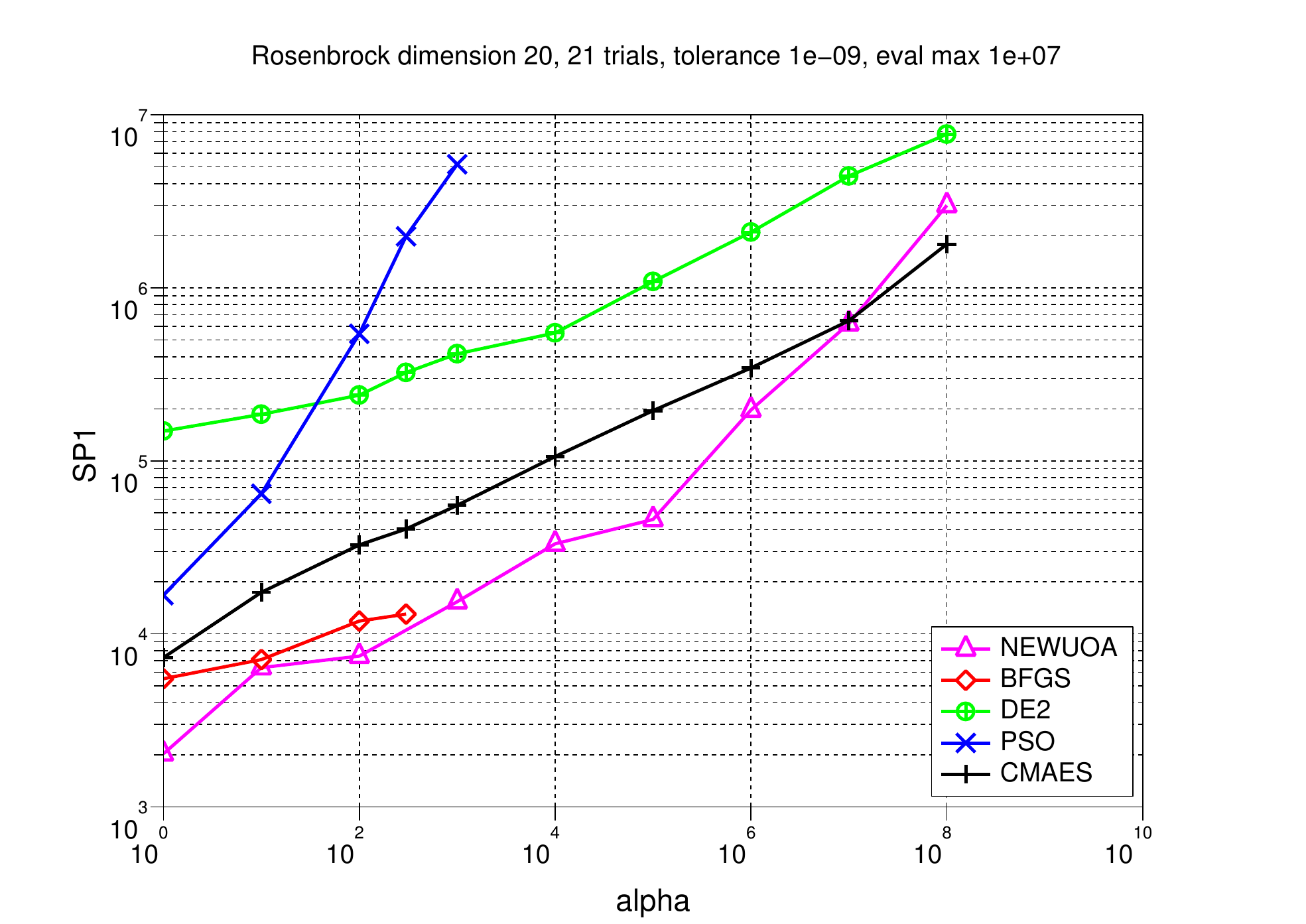}
\end{minipage}
\begin{minipage}{0.49\textwidth}\centering
Rotated Rosenbrock Function\\
\includegraphtrim{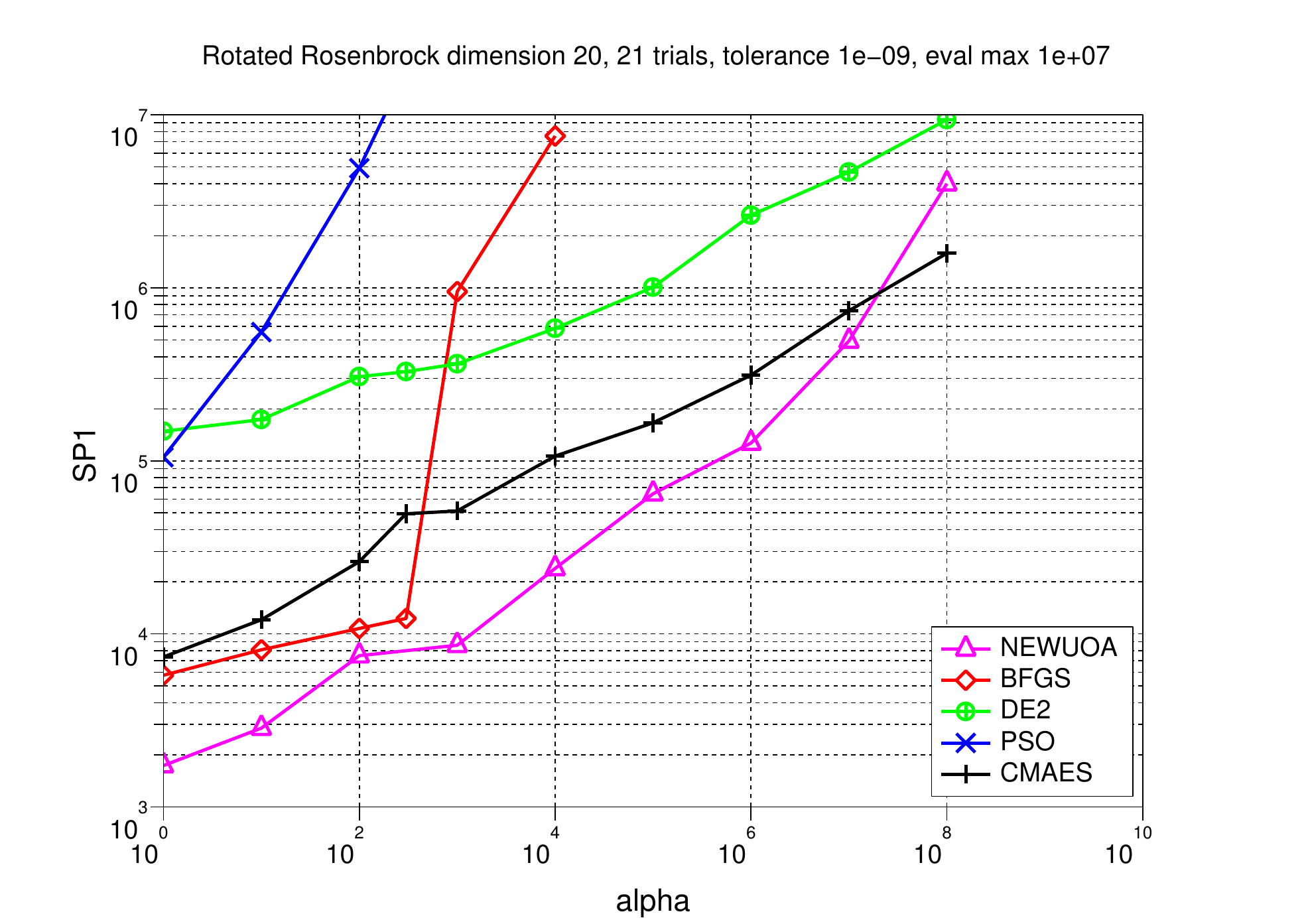}
\end{minipage}
\caption[rosen]{\label{fig:rosen}Rosenbrock function. Shown is $SP1$ (the expected running time or number of function evaluations to reach the target function value) versus conditioning parameter $\alpha$. }
\end{figure}
%
On the Rosenbrock function \NEW\ is the best algorithm (Figure~\ref{fig:rosen}). \NEW\ outperforms CMA-ES roughly by a factor of five, vanishing for very large values for the conditioning parameter \alp.  For small \alp, BFGS is in-between, and for $\alp>10^4$ BFGS fails. DE is again roughly ten times slower than CMA-ES. Only PSO shows a strong dependency on the rotation of the function and it reveals the strongest performance decline with increasing \alp, while it never competes with the best three algorithms.

\paragraph{Scaling behaviors} The scaling of the performance with search space dimension is similar for all functions (see Appendix~\ref{app:results} for the data). CMA-ES, \NEW\ and PSO show the best scaling behavior. They slow down by a factor between five and ten in 40D compared to 10D. For BFGS the factor is slightly above ten, while for DE the factor is thirty or larger, presumably because the default population size increase linearly with the dimension.

\section{Summary}
\label{conclusions}
In this paper we have conducted a comparison of BFGS, NEWUOA, and three stochastic bio-inspired optimization methods in a black-box optimization scenario. The empirical study was conducted on smooth functions with varying condition number. Aside from gradients being not provided, we consider these functions as the favorite playgrounds of BFGS and NEWUOA. We find that NEWUOA performs exceptional on separable quadratic functions, it performs in all cases very well with moderate condition numbers, but shows a comparatively steep performance decline with increasing ill-conditioning. BFGS performs well overall, but shows a strong decline on very ill-conditioned non-separable functions. For DE, the parameters are difficult to tune and yet it performs overall poorly with the single best parameter setting on our small function set. With the chosen parameters, DE shows the strongest robustness to ill-conditioning though. PSO performs similar to CMA-ES on the separable problems, with an even weaker dependency on the conditioning. On non-separable problems PSO exhibits a strong performance decline with increasing conditioning and performs very poorly even on moderately ill-conditioned functions.  Finally, CMA-ES generally outperforms DE and PSO, while up to a moderate function conditioning BFGS and NEWUOA are significantly faster in most cases. Due to their invariance properties, the performance results of CMA-ES and DE are the most stable ones and most likely to generalize to other functions. 

\section*{Acknowledgements}
We would like to acknowledge Philippe Toint for his kind suggestions, and Nikolas Mauny for writing the Scilab transcription of the Standard PSO 2006 code.

\bibliography{bibliography}
\newpage
\begin{appendix}
\section{All Results}\label{app:results}
\begin{figure}[ht]
\begin{minipage}{0.49\textwidth}\centering
Separable Ellipsoid Function\\
\includegraphtrim{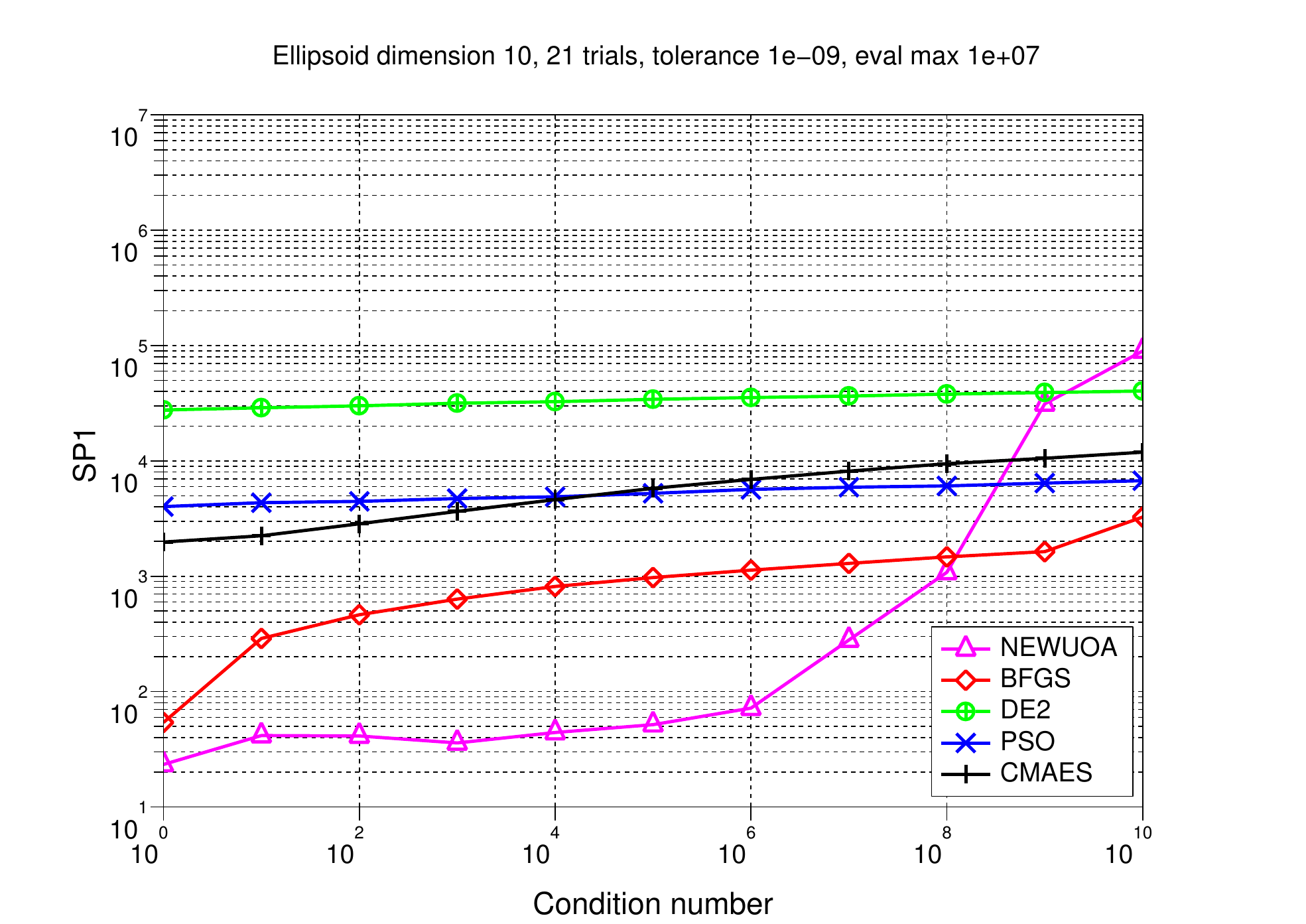}
\includegraphtrim{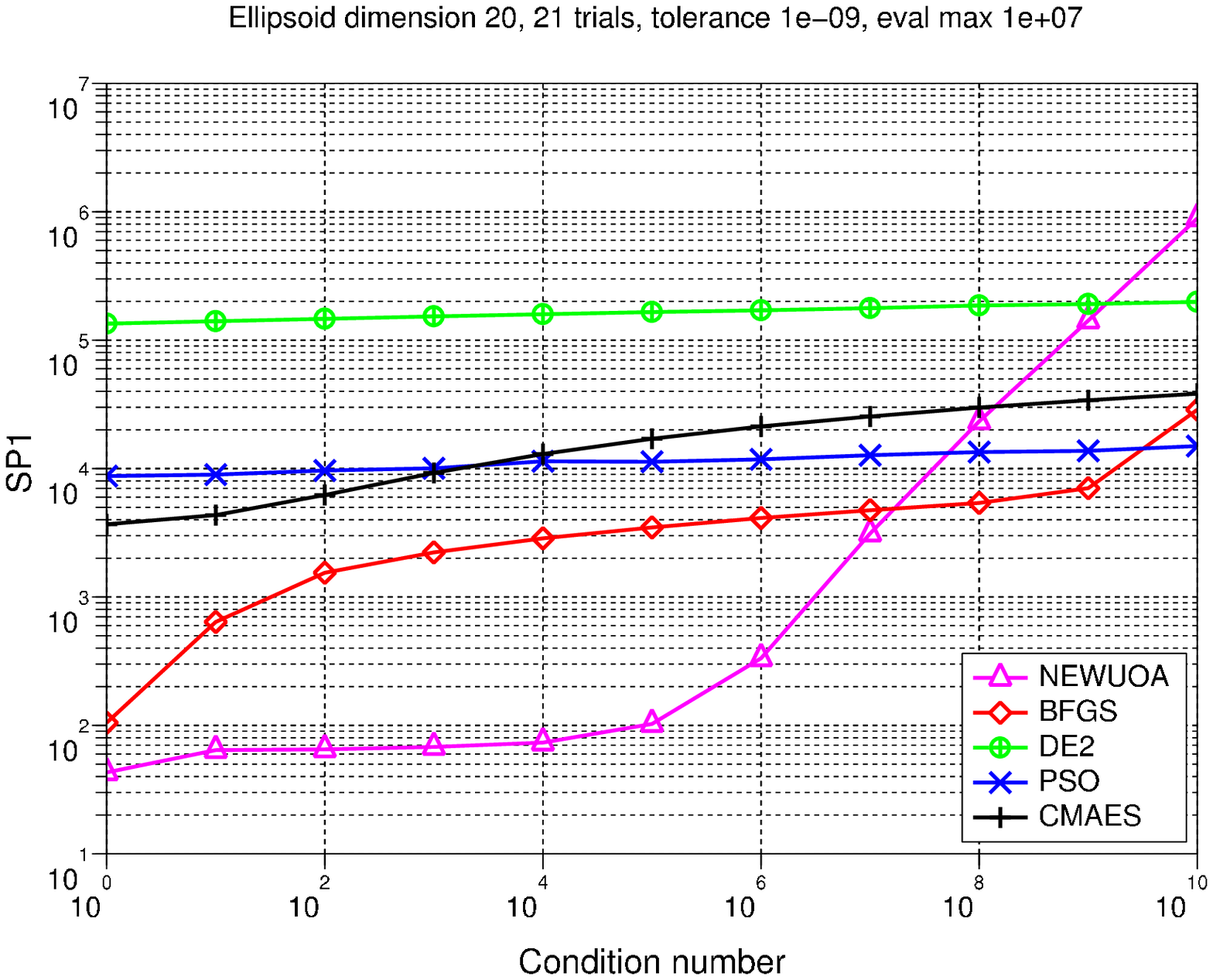}
\includegraphtrim{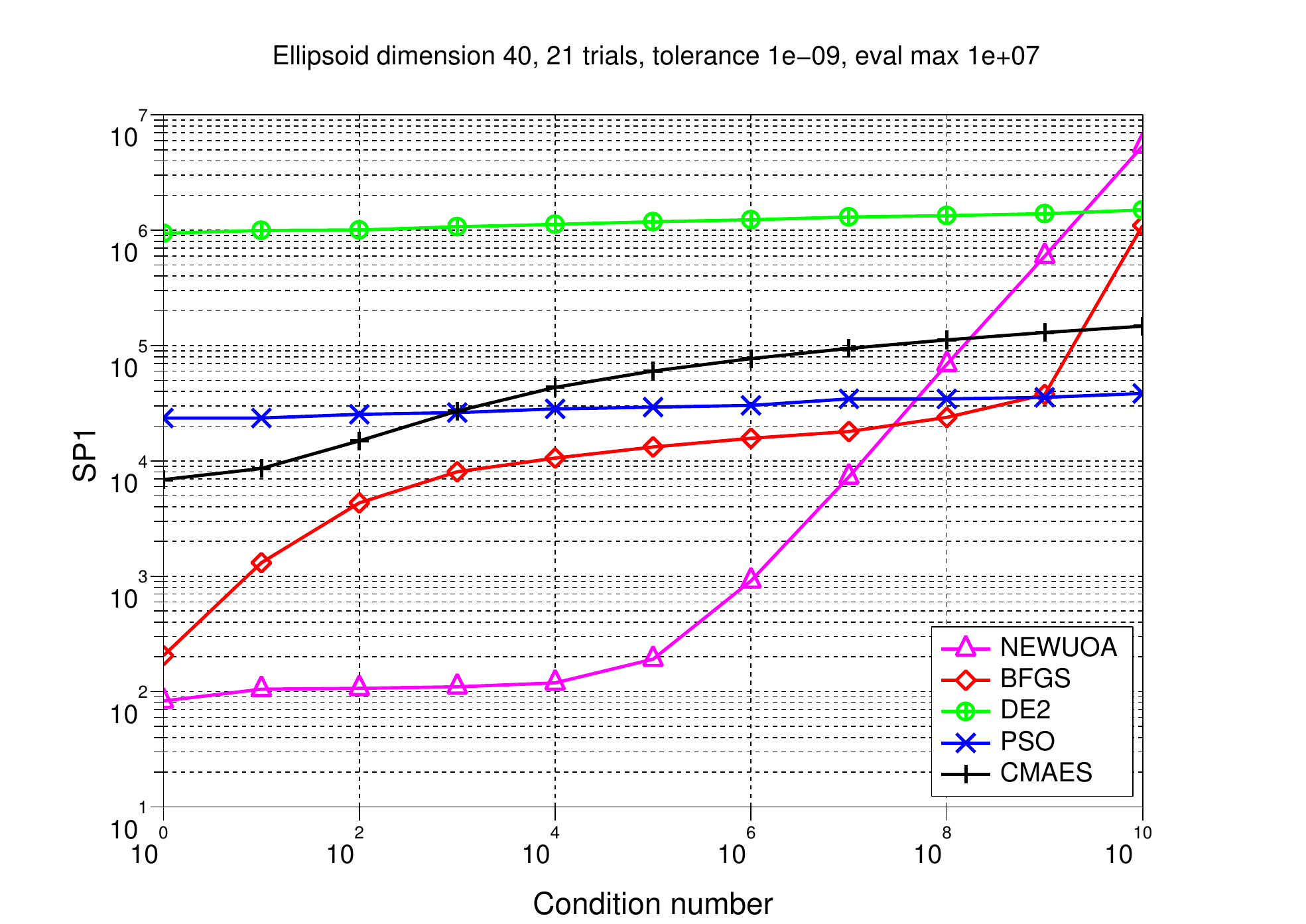}
\end{minipage}
\begin{minipage}{0.49\textwidth}\centering
Rotated Ellipsoid Function\\
\includegraphtrim{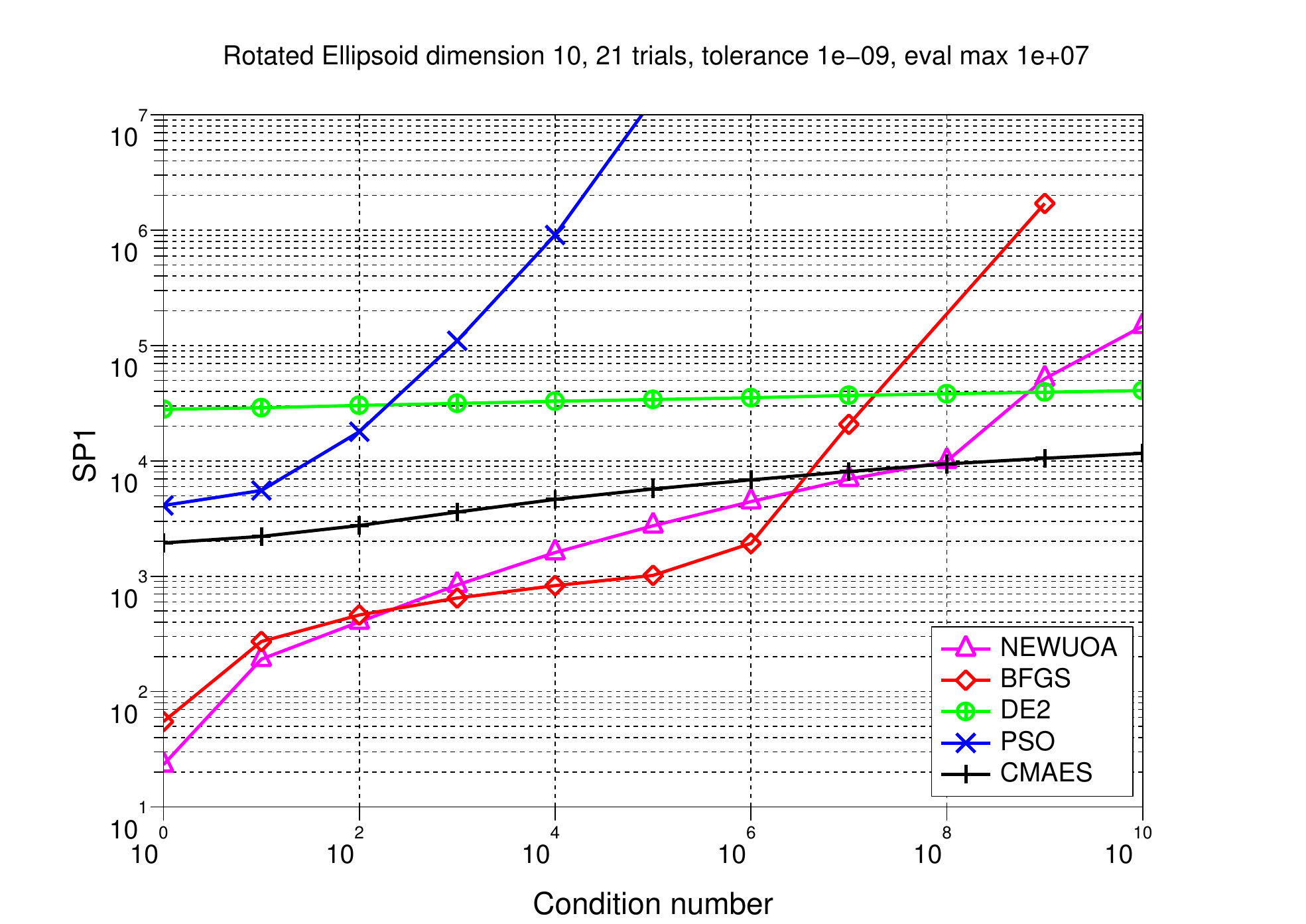}
\includegraphtrim{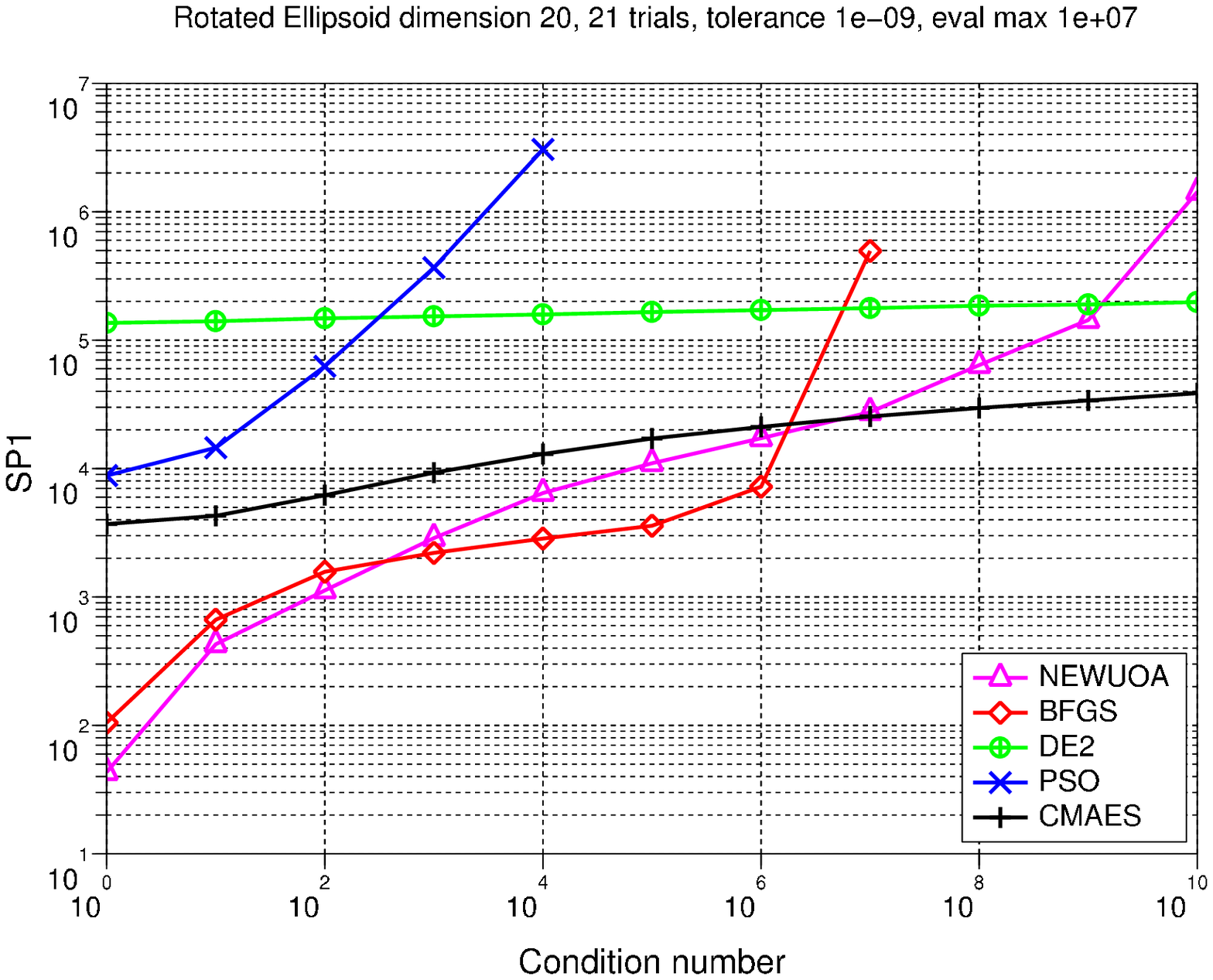}
\includegraphtrim{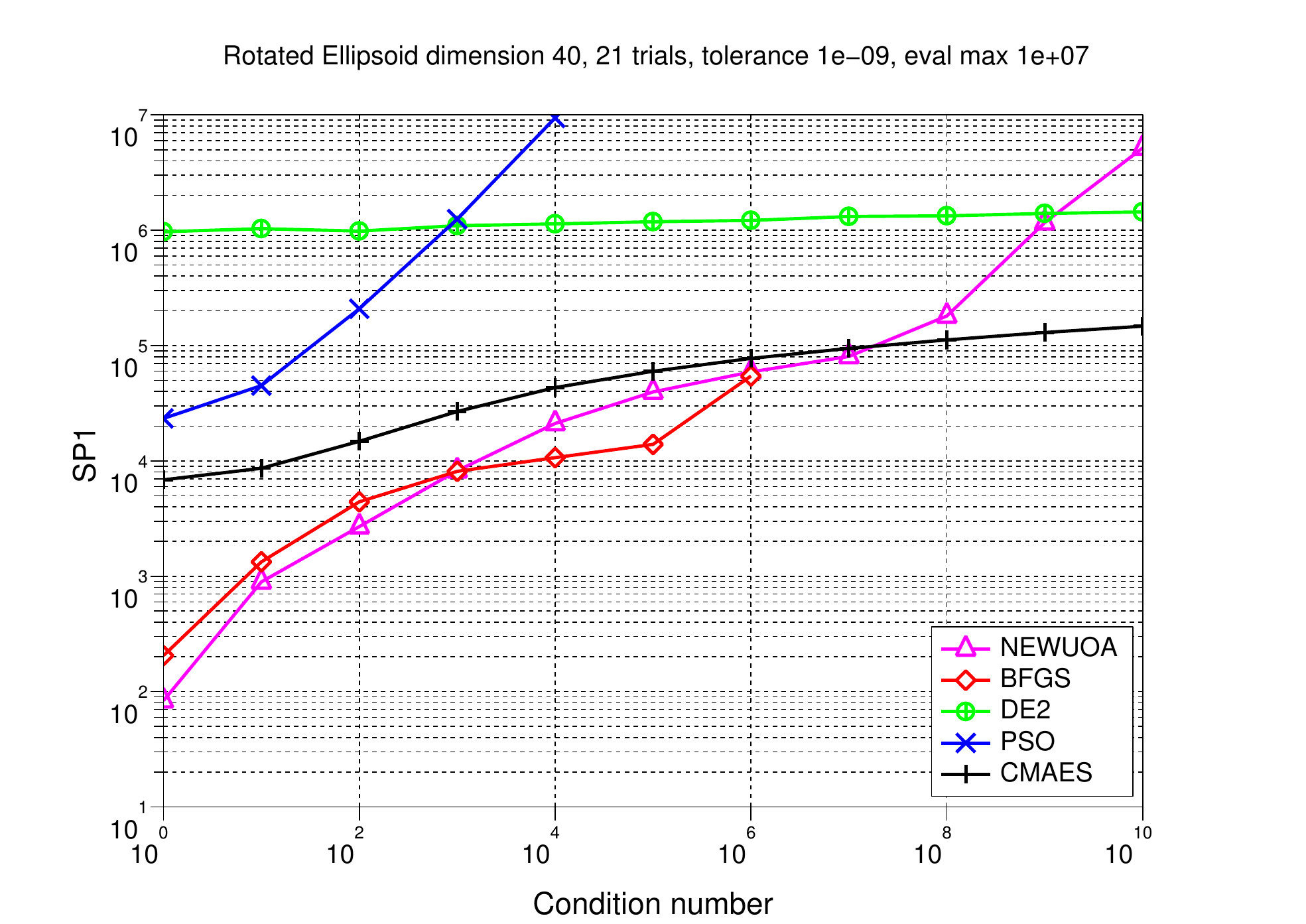}
\end{minipage}
 \caption[elli]{Ellipsoid function. Shown is $SP1$ (the expected running time or number of function evaluations to reach the target function value) versus condition number.}
\end{figure}

\begin{figure}[ht]
\begin{minipage}{0.49\textwidth}\centering
Separable \SElli\ Function\\
\includegraphtrim{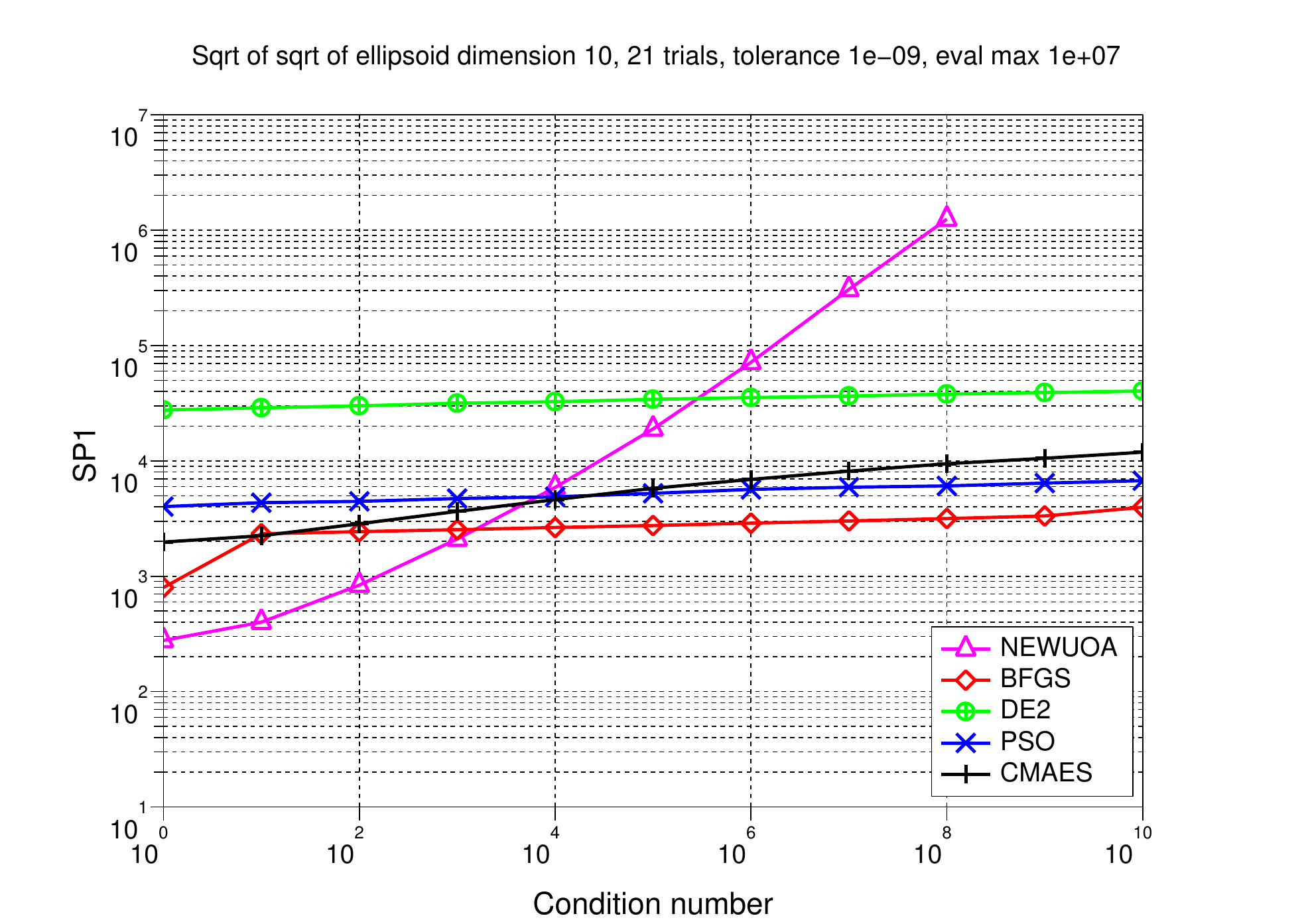}
\includegraphtrim{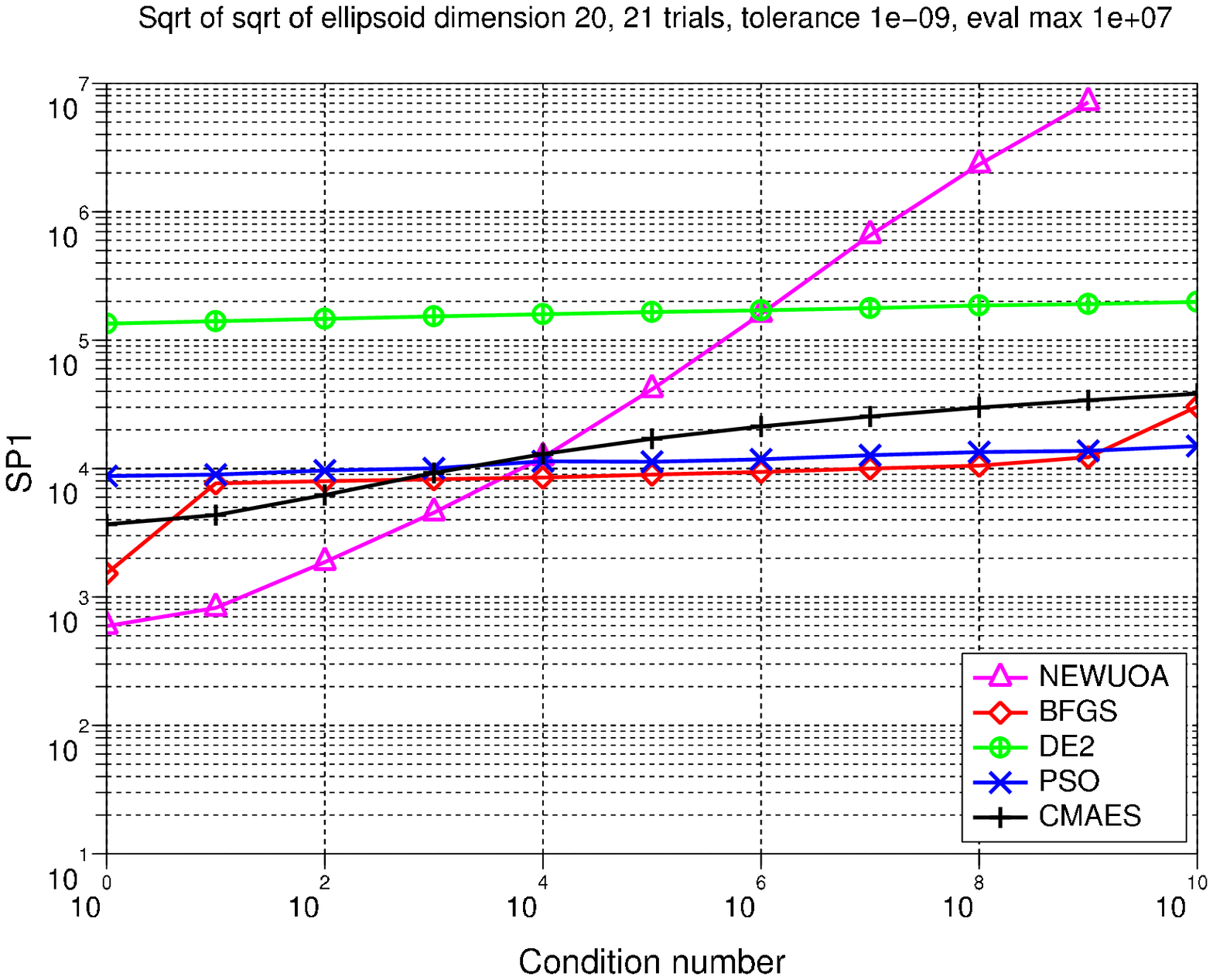}
\includegraphtrim{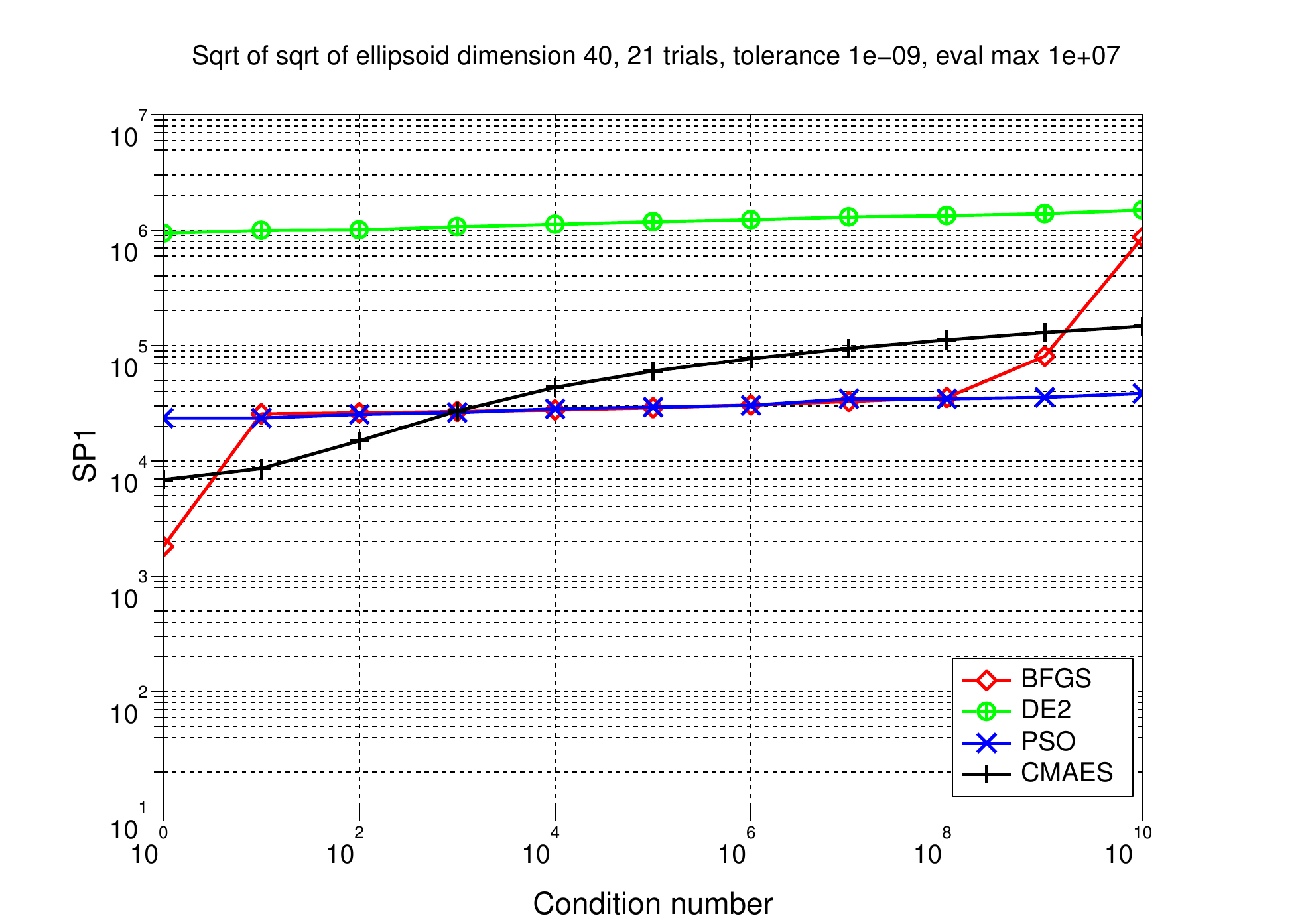}
\end{minipage}
\begin{minipage}{0.49\textwidth}\centering
Rotated \SElli\ Function\\
\includegraphtrim{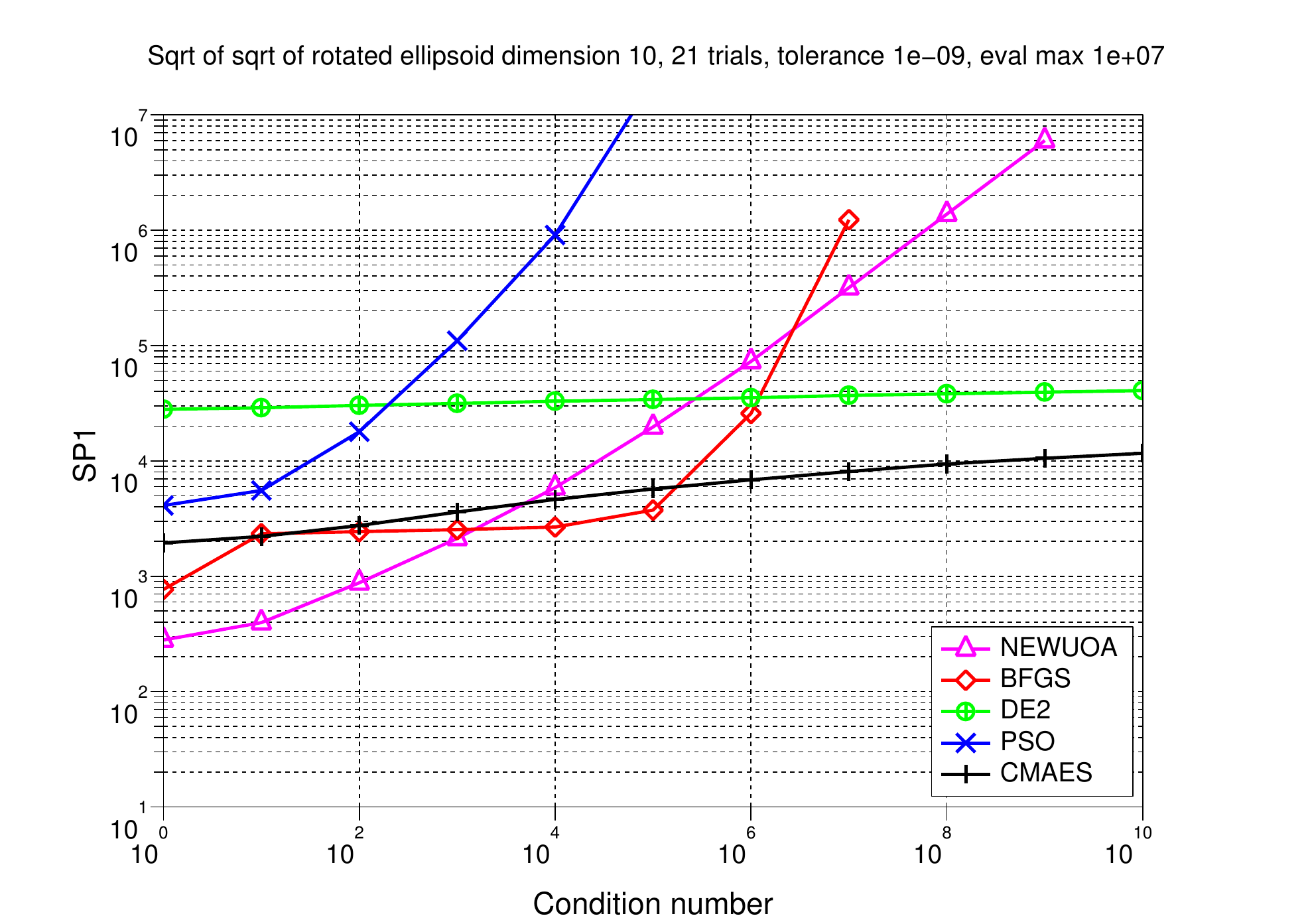}
\includegraphtrim{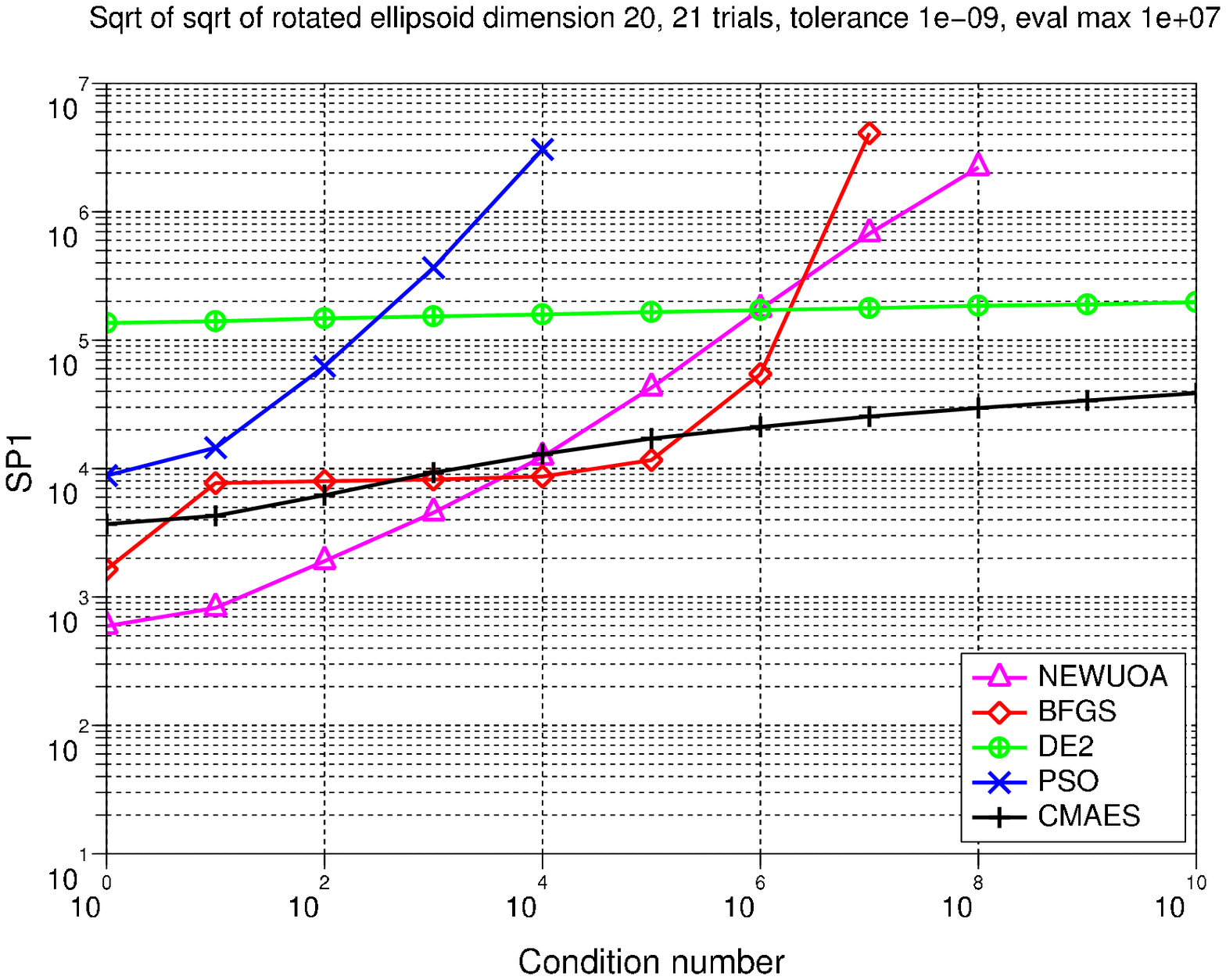}
\includegraphtrim{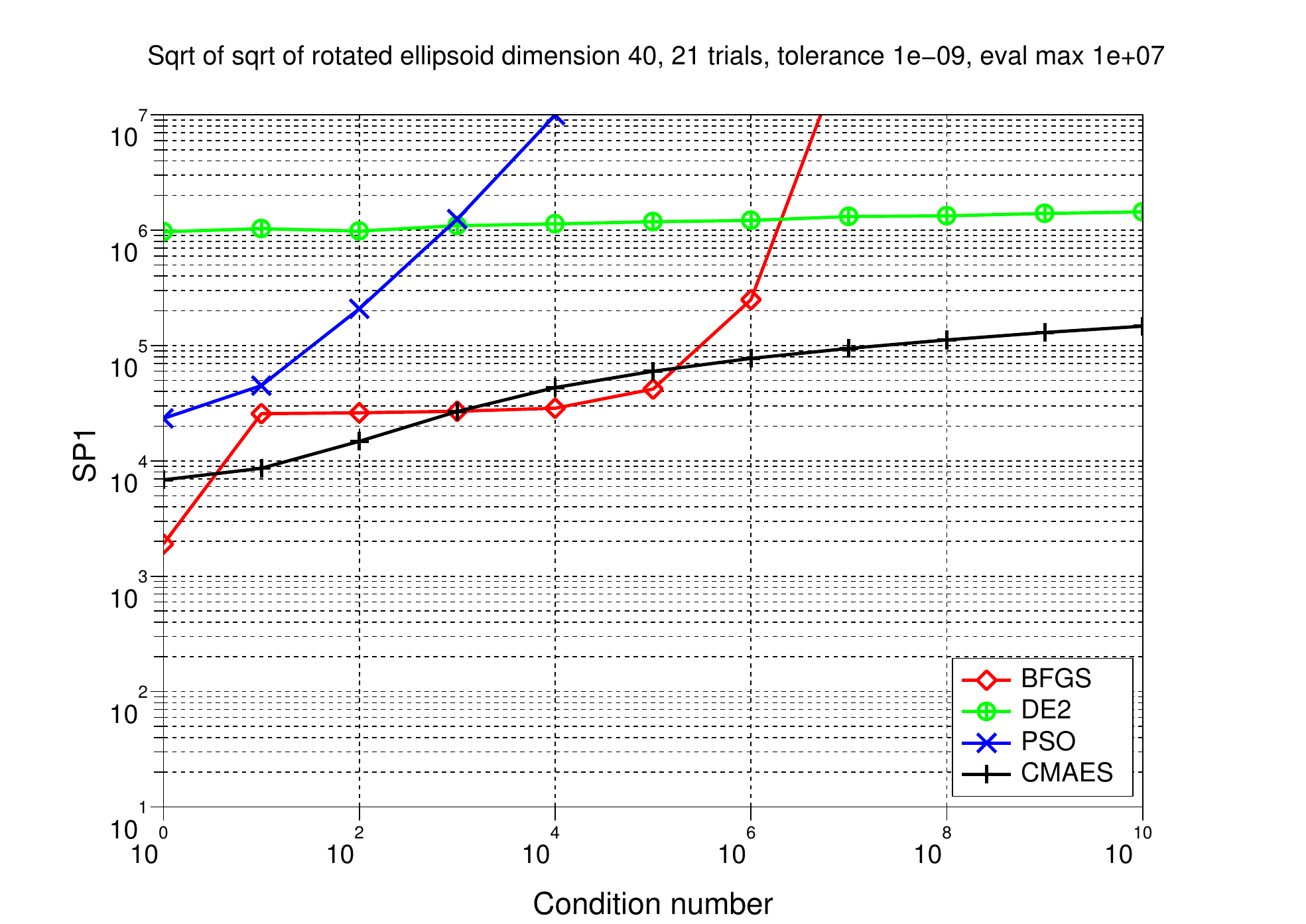}
\end{minipage}
\caption[elli]{Ellipsoid$^{1/4}$ function. Shown is $SP1$ (the expected running time or number of function evaluations to reach the target function value) versus condition number.}
\end{figure}

\begin{figure}[t]
\begin{minipage}{0.49\textwidth}\centering
Rosenbrock Function\\
\includegraphtrim{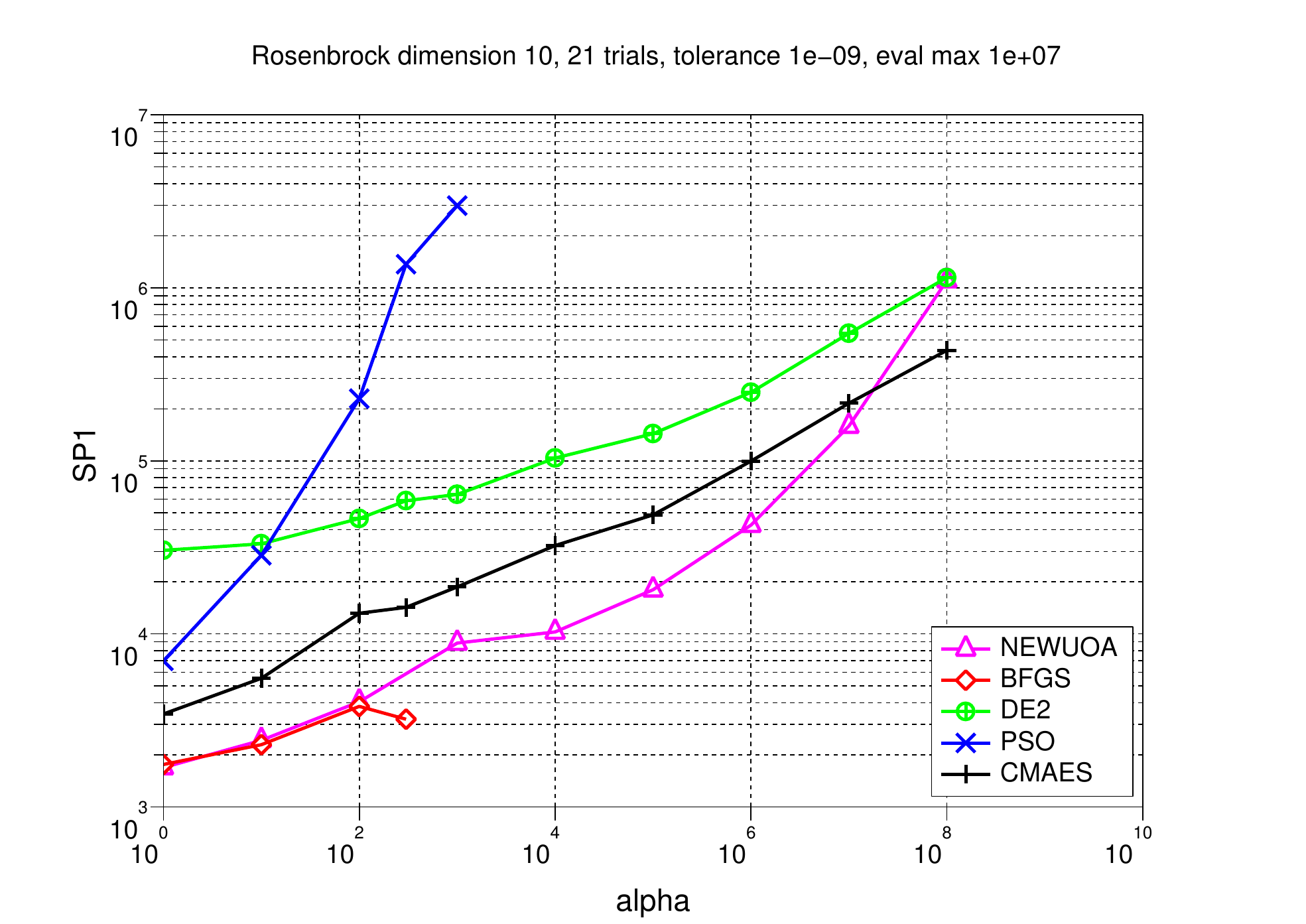}
\includegraphtrim{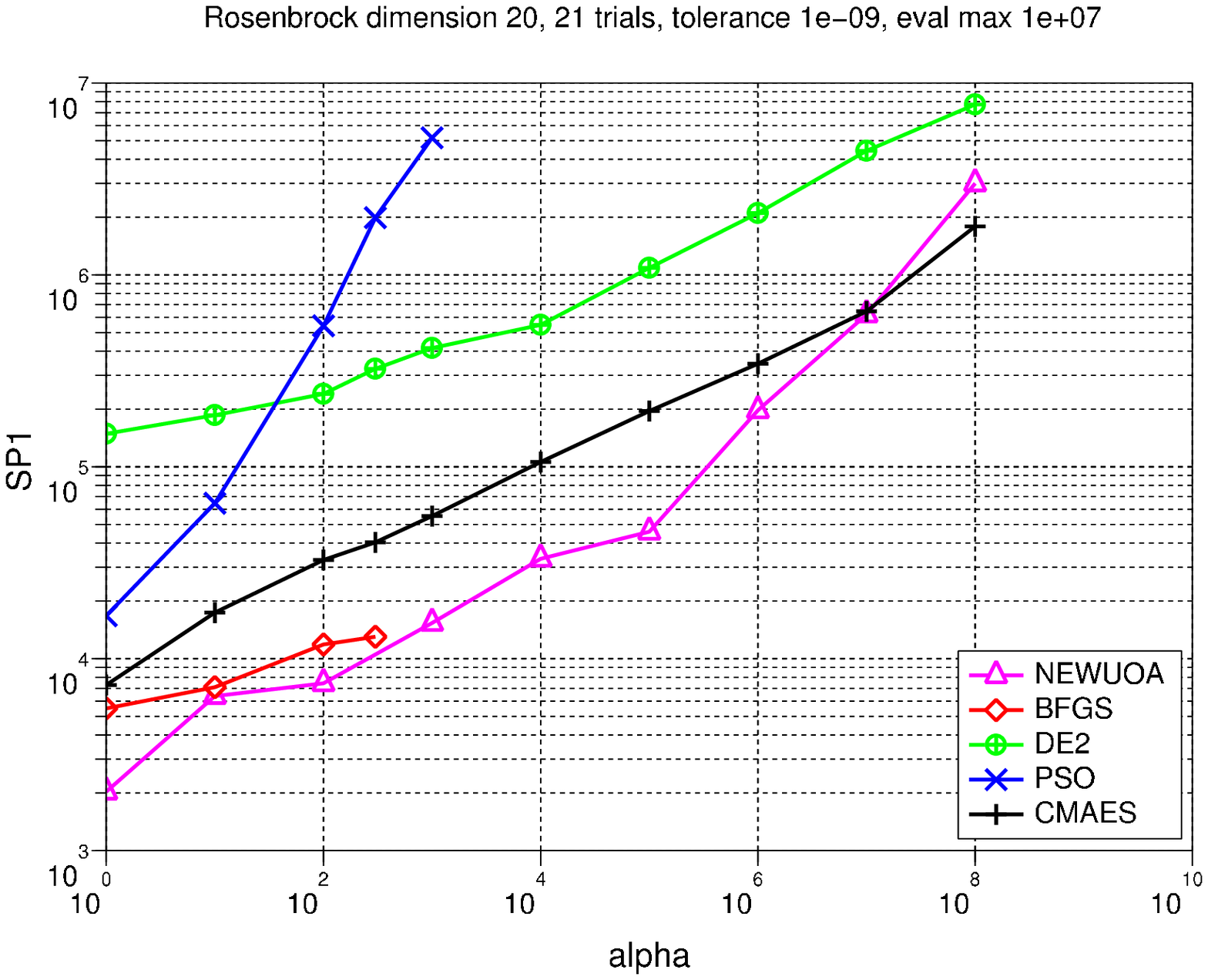}
\includegraphtrim{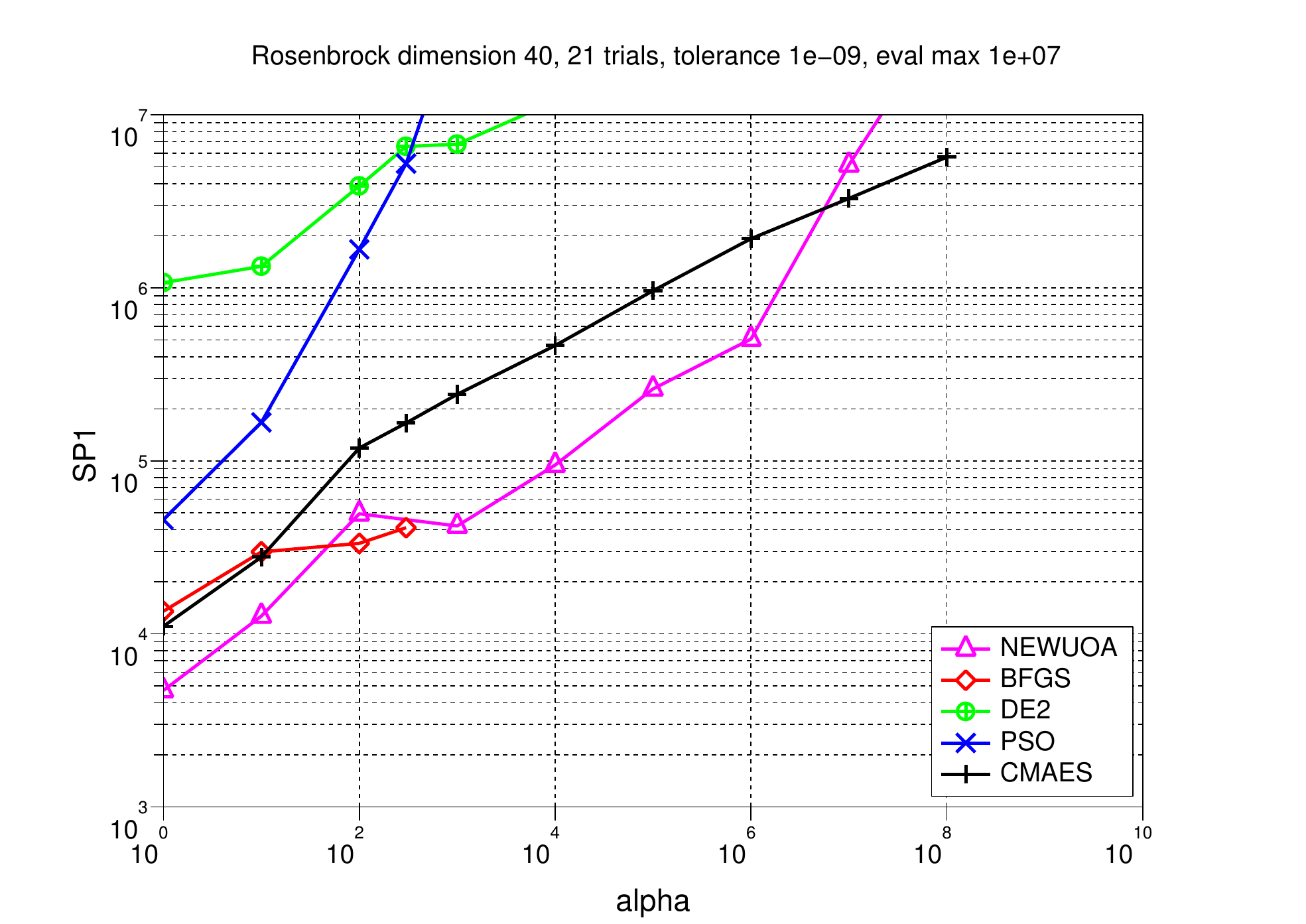}
\end{minipage}
\begin{minipage}{0.49\textwidth}\centering
Rotated Rosenbrock Function\\
\includegraphtrim{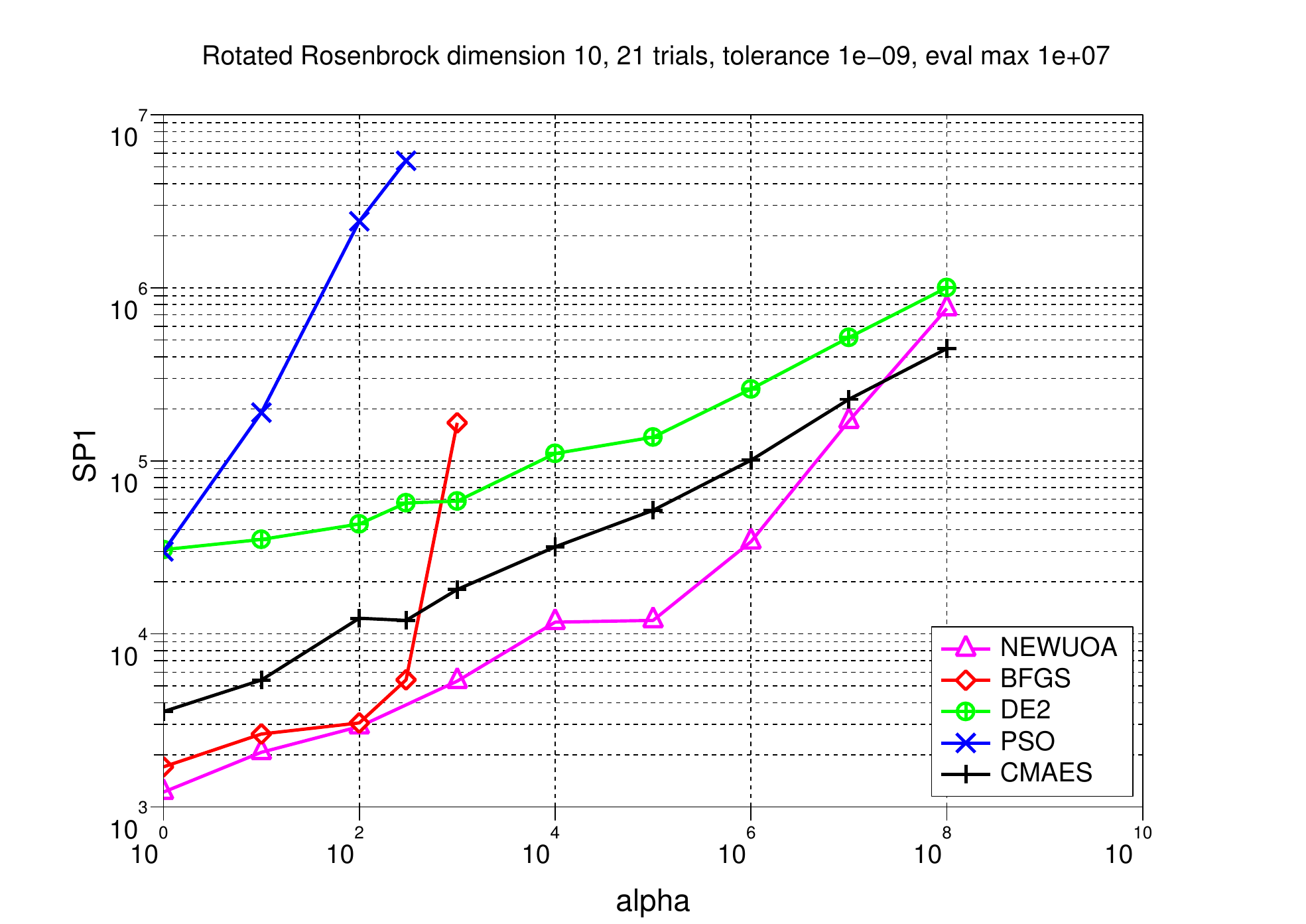}
\includegraphtrim{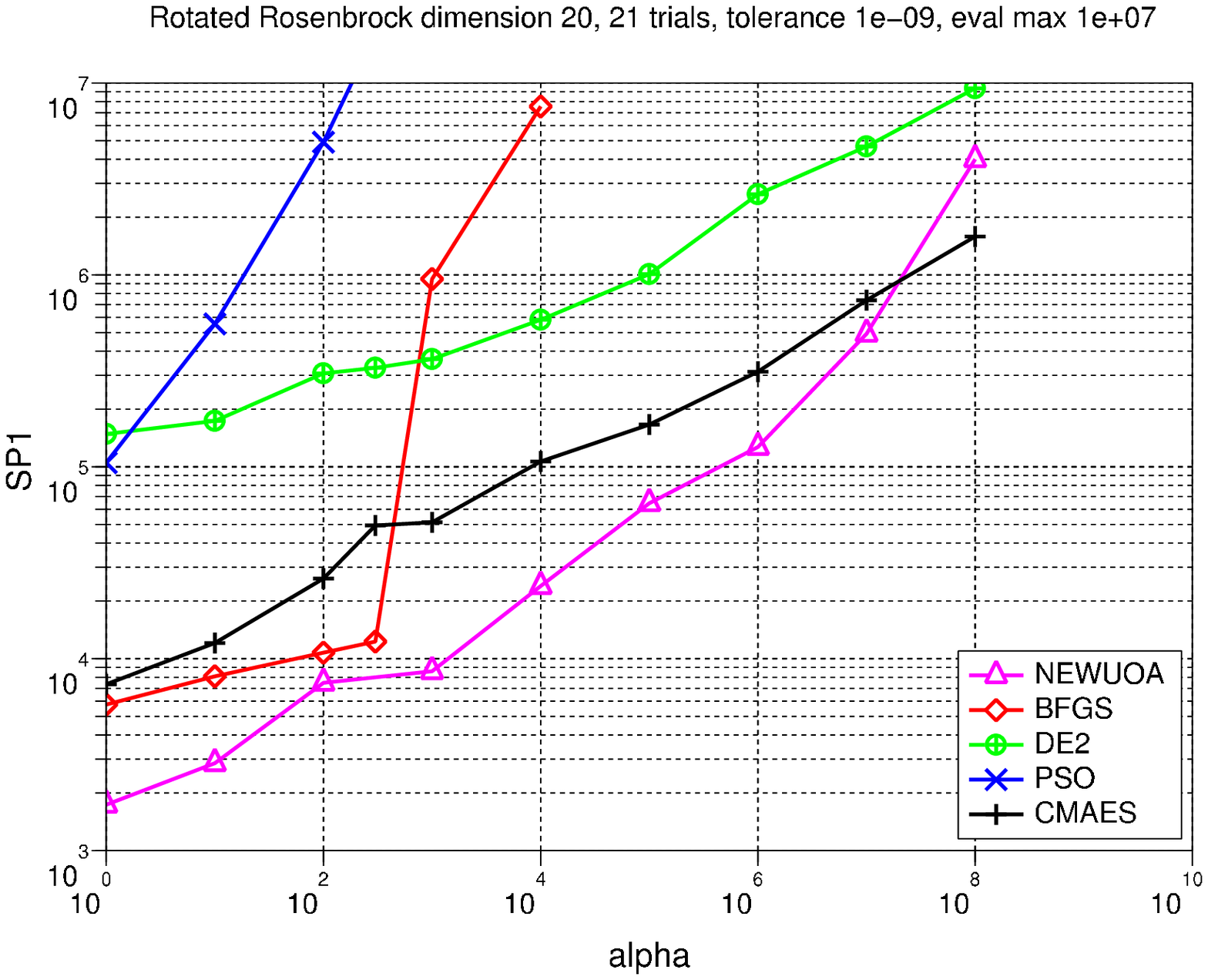}
\includegraphtrim{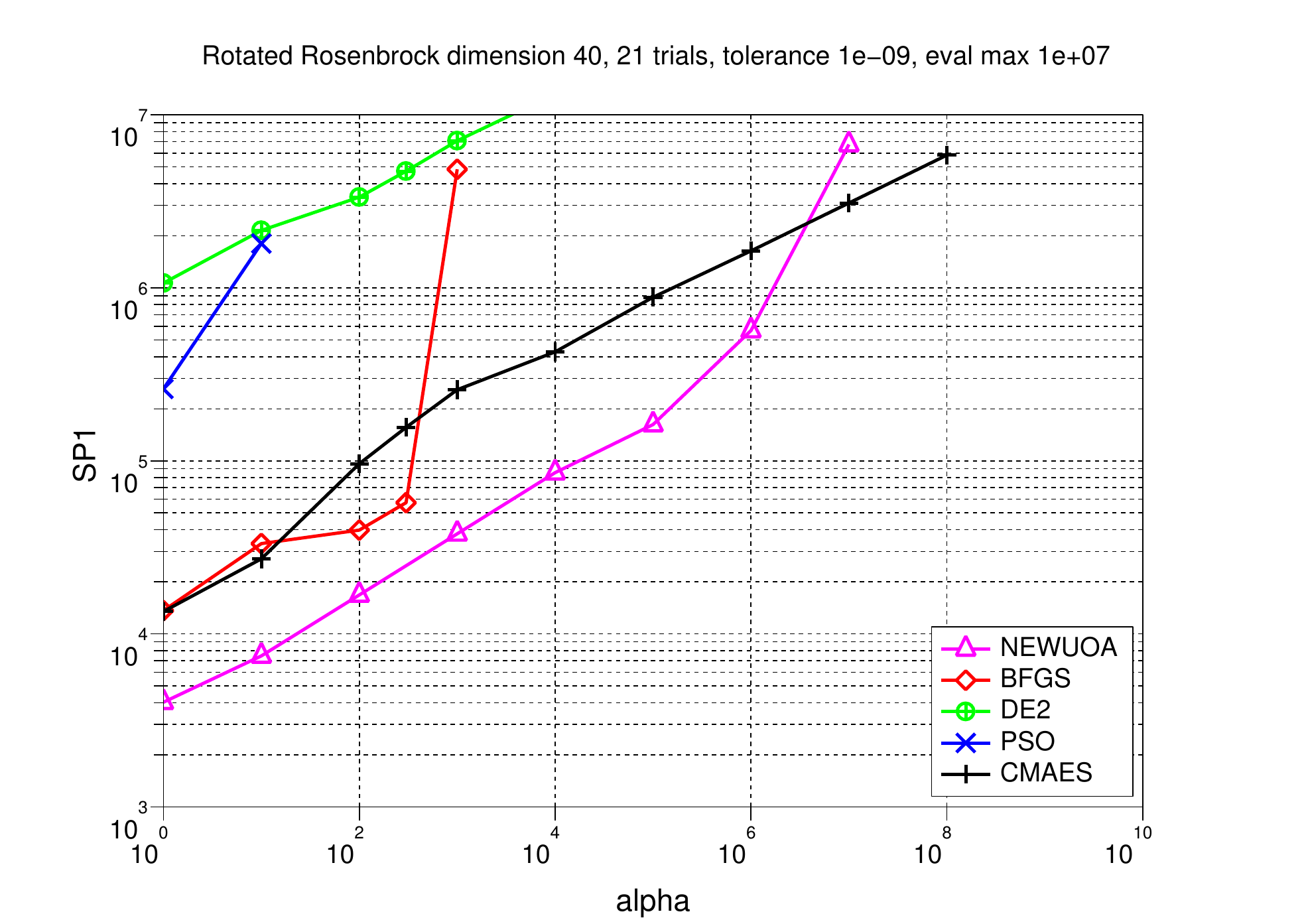}
\end{minipage}
\caption[rosen]{Rosenbrock function. Shown is $SP1$ (the expected running time or number of function evaluations to reach the target function value) versus conditioning parameter $\alpha$. }
\end{figure}

\end{appendix}
\end{document}